\documentclass[11pt]{article}
\usepackage{amssymb}
\usepackage{color,geometry}

\newcommand{\R}{\mathbb{R}}

\newcommand{\cuad}{{\sqcap\kern-.68em\sqcup}}
\newcommand{\abs}[1]{|#1|}
\newcommand{\norm}[1]{\|#1\|}

\newtheorem{definition}{Definition}[section]

\newtheorem{proposition}{Proposition}[section]

\newtheorem{lemma}{Lemma}[section]

\newtheorem{remark}{Remark}[section]
\newcommand{\bremark}{\begin{remark} \em}
\newcommand{\eremark}{\end{remark} }

\newcommand{\be}{\begin{equation}}
\newcommand{\bel}[1]{\begin{equation}\label{#1}}
\newcommand{\ee}{\end{equation}}
\newcommand{\barr}{\begin{eqnarray}}
\newcommand{\earr}{\end{eqnarray}}
\newcommand{\bars}{\begin{eqnarray*}}
\newcommand{\ears}{\end{eqnarray*}}

\newtheorem{subn}{\name}

\newcommand{\bsn}[1]{\def\name{#1}\begin{subn}}
\newcommand{\esn}{\end{subn}}
\newtheorem{sub}{\name}[section]
\newcommand{\bs}{\begin{sub}}
\newcommand{\es}{\end{sub}}
\newcommand{\bsl}[1]{\begin{sub}\label{#1}}
\newcommand {\dint}{{\displaystyle \int\!\!\!\int}}
\newcommand{\bth}[1]{\def\name{Theorem}
\begin{sub}\label{t:#1}}
\newcommand{\blemma}[1]{\def\name{Lemma}
\begin{sub}\label{l:#1}}
\newcommand{\bcor}[1]{\def\name{Corollary}
\begin{sub}\label{c:#1}}
\newcommand{\bdef}[1]{\def\name{Definition}
\begin{sub}\label{d:#1}}
\newcommand{\bprop}[1]{\def\name{Proposition}
\begin{sub}\label{p:#1}}
\newcommand{\rth}[1]{Theorem~\ref{t:#1}}
\newcommand{\rlemma}[1]{Lemma~\ref{l:#1}}
\newcommand{\rcor}[1]{Corollary~\ref{c:#1}}

\newcommand{\rprop}[1]{Proposition~\ref{p:#1}}
\newcommand{\BA}{\begin{array}}
\newcommand{\EA}{\end{array}}
\newcommand{\BAN}{\renewcommand{\arraystretch}{1.2}
\setlength{\arraycolsep}{2pt}\begin{array}}
\newcommand{\BAV}[2]{\renewcommand{\arraystretch}{#1}
\setlength{\arraycolsep}{#2}\begin{array}}
\newcommand{\BSA}{\begin{subarray}}
\newcommand{\ESA}{\end{subarray}}
\newcommand{\BAL}{\begin{aligned}}
\newcommand{\EAL}{\end{aligned}}
\newcommand{\BALG}{\begin{alignat}}
\newcommand{\EALG}{\end{alignat}}
\newcommand{\BALGN}{\begin{alignat*}}
\newcommand{\EALGN}{\end{alignat*}}
\newcommand{\note}[1]{\textit{#1.}\hspace{2mm}}
\newcommand{\Proof}{{\bf Proof.\ }}
\newcommand{\qeda}{\hspace{10mm}\hfill $\square$}

\newcommand{\Remark}{\note{Remark}}

\def\angb<#1>{\langle #1 \rangle}

\newcommand{\opname}[1]{\mbox{\rm #1}\,}

\newcommand{\dist}{\opname{dist}}
\newcommand{\myfrac}[2]{{\displaystyle \frac{#1}{#2} }}
\newcommand{\myint}[2]{{\displaystyle \int_{#1}^{#2}}}

\newcommand{\prt}{\partial}

\newcommand{\ti}{\times}

\newcommand{\nind}{\noindent}

\def\ga{\alpha}     \def\gb{\beta}       \def\gg{\gamma}
       \def\gd{\delta}      \def\ge{\epsilon}
\def\gth{\theta}                         
\def\gf{\phi}           
      \def\gk{\kappa}      \def\gl{\lambda}
\def\gm{\mu}        \def\gn{\nu}         \def\gp{\pi}
    \def\gr{\rho}        
\def\gs{\sigma}       \def\gt{\tau}
\def\gu{\upsilon}      \def\gw{\omega}
                \def\gz{\zeta}
\def\Gg{\Gamma}     \def\Gd{\Delta}      \def\Gf{\Phi}
\def\Gth{\Theta}
\def\Gl{\Lambda}    \def\Gs{\Sigma}

\def\CS{{\mathcal S}}   \def\CM^+{{\mathcal M}}   
\def\CR{{\mathcal R}}      
      
      \def\CF{{\mathcal F}}
      
      \def\CL{{\mathcal L}}

   \def\BBE {\mathbb E}    
   \def\BBH {\mathbb H}    
       
   \def\BBN {\mathbb N}    
   \def\BBR {\mathbb R}    
       
       \def\BBY {\mathbb Y}

\begin{document}

\begin{center}{\bf  \Large    Initial trace of positive solutions to fractional \\[2mm] diffusion equations with absorption }
\bigskip
\medskip

{\small

 {\bf Huyuan Chen}\footnote{chenhuyuan@yeah.net}   \qquad     {\bf Laurent V\'eron}\footnote{laurent.veron@lmpt.univ-tours.fr}\medskip \smallskip

     $^1$Department of Mathematics, Jiangxi Normal University,
       \\ Nanchang 330022, China \\[2mm]

       \smallskip
$^2$Laboratoire de Math\'{e}matiques et Physique Th\'{e}orique
\\  Universit\'{e} de Tours, 37200 Tours, France\\[1mm]

}
\medskip

\begin{abstract}
In this paper, we prove the existence of an initial trace ${\cal T}_u$ for any positive solution $u$ to the semilinear fractional diffusion equation $(H)$
$$\partial_t u + (-\Delta)^s  u+f(t,x,u)=0\quad  {\rm in}\quad (0,+\infty)\times\R^N,$$
where $N\geq1$,  the operator $(-\Delta)^s$ with $s\in(0,1)$ is the fractional Laplacian, $f:\R_+\times\R^N\times\R_+\rightarrow\BBR$ is a Caratheodory function satisfying $f(t,x,u)u\geq 0$ for all $(t,x,u)\in \R_+\times\R^N\times\R_+$ and $\R_+=[0,+\infty)$.
We define the regular set of the trace ${\cal T}_u$ as an open subset of ${\cal R}_u\subset\BBR^N$ carrying a nonnegative Radon measure $\nu_u$ such that
$$\lim_{t\to 0}\myint{{\cal R}_u}{}u(t,x)\gz(x)\, dx=\myint{{\cal R}_u}{}\gz\, d\gn_u,\quad\forall\, \gz\in C^{2}_0({\cal R}_u),
$$
and the singular set ${\cal S}_u=\R^N\setminus {\cal R}_u$ as the set points $a$ such that
$$\limsup_{t\to 0}\myint{B_\gr(a)}{}u(t,x)\, dx=+\infty\quad{\rm for\ any}\ \, \gr>0.
$$
We also study the reverse problem of constructing a positive solution to $(H)$ with a given initial trace $(\CS,\gn)$, where $\CS\subset\BBR^N$ is a closed set and $\gn$ is a positive Radon measure on $\CR=\BBR^N\setminus\CS$ and develop the case  $f(t,x,u)=t^\gb u^p$ with $\gb>-1$ and
$p>1$.

\end{abstract}
\end{center}


  \vspace{1mm}


\tableofcontents
\hspace{.05in}

\vspace{2mm}


\setcounter{equation}{0}
\section{Introduction}

The first aim of this paper is to study the existence of an {\it initial trace} of positive solutions  to the semilinear fractional diffusion equation
\begin{equation}\label{In-1}
\partial_t u + (-\Delta)^s  u+f(t,x,u)=0\quad  {\rm in}\quad Q_\infty:=\BBR^*_+\times\R^N,
\end{equation}
where   $f:\R^*_+\times\R^N\times\R\rightarrow\BBR$ is a Caratheodory function satisfying
\bel{In-2}
f(t,x,u)u\geq 0,\quad\forall\,(t,x,u)\in \R^*_+\times\R^N\times\R,
\ee
and  $\R^*_+=(0,+\infty)$.
The fractional Laplacian$(-\Delta)^s $ with $s\in(0,1)$ is defined in the principal value sense that
$$(-\Delta)^s  u(x)=\lim_{\varepsilon\to0^+} (-\Delta)_\varepsilon^s u(x),$$
where
\begin{equation}\label{In-0}
(-\Delta)_\varepsilon^s  u(x):=-a_{N,s}\int_{\R^N}\frac{ u(z)-
u(x)}{|z-x|^{N+2s}}\chi_\varepsilon(|x-z|) dz\;,\;\; a_{N,s}=\myfrac{\Gg(\frac N2+s)}{\gp^{\frac N2}\Gg(2-s)}s(1-s),
\end{equation}
for $\varepsilon>0$ and
$$\chi_\varepsilon(r)=\left\{ \arraycolsep=1pt
\begin{array}{lll}
0 \quad & {\rm  if}\quad r\in[0,\varepsilon]\\[1mm]
1 \quad & {\rm  if}\quad r>\varepsilon.
\end{array}
\right.$$
The solutions of (\ref{In-1}) are intended in the classical sense and, in order $(-\Gd)^{s}u(t,x)$ to be well-defined, we always assume that $u(t,.)\in \CL^s(\BBR^N)$ for any $t>0$, where
\bel{In-X}
\CL^s(\BBR^N)=\left\{\phi\in L^1_{loc}(\BBR^N)\;\;{\rm s.t. }\;\;\norm\gf_{\CL^s}:=\myint{\BBR^N}{}\myfrac{\abs{\gf(x)}\,dx}{1+\abs{x}^{N+2s}}<+\infty\right\}.
\ee
Notice that the constant functions belong to $\CL^s(\BBR^N)$. If $\gw\subset\BBR^N$ and $0<T\leq+\infty$, we set $Q_T^\gw=(0,T)\ti\gw$, $Q_T^{\BBR^N}=Q_T$, $Q_\infty=\BBR_+^*\times\BBR^N$ and denote by $B_\gr(z)$ (resp. $K_\gr(z)$) the open ball (resp. open cube with sides parallel to the axis) with center $z\in\BBR^N$ and radius (side length) $\gr>0$. We define the {\it regular set} of the initial trace of a positive solution $u$ of (\ref{In-1}) by
\bel{In-3}
\CR_u=\left\{z\in\BBR^N:\exists \,\gr>0\;\,{\rm  s.t. }\, \dint_{\!Q_1^{B_\gr(z)}}f(t,x,u)\,dxdt<+\infty\right\}.
\ee
Clearly $\CR_u$ is open. The {\it conditional singular set} $\tilde \CS_u$ is $\BBR^N\setminus\CR_u$ and {\it  the conditional initial trace} is the couple $Tr_c(u):=(\tilde \CS_u,\gn)$. Our first result is the following statement which is the starting point of our work.\medskip

\nind{\bf Theorem A} {\it Let $u$ be a nonnegative classical solution of (\ref{In-1}) and
the regular set $\CR_u$ of $u$ is given in (\ref{In-3}), then  there exists a nonnegative Radon measure $\gn_u$ on $\CR_u$ such that
\bel{In-4}
\lim_{t\to 0}\myint{\CR_u}{}u(t,x)\gz(x)\, dx=\myint{\CR_u}{}\gz \,d\gn_u,\quad\forall\, \gz\in C^2_0(\CR_u).
\ee
}\medskip

The problem of the initial trace of nonnegative solutions for semilinear heat equations was initiated by Marcus and V\'eron in \cite{MV1} with equation
\bel{In-6}
\prt_tu-\Gd u+u^p=0\quad{\rm in }\;\;Q_\infty,
\ee
for $p>1$. They showed the existence of an initial trace $Tr(u)$ represented by a closed subset $\CS_u$ of $\BBR^N$ and a nonnegative Radon measure $\gn_u$ on $\CR_u=\BBR^N\setminus\CS_u$. On $\CR_u$ the initial trace is achieved as in $(\ref{In-4})$. On $\CS_u$ they proved that for any $z\in\CS_u$,
\bel{In-7}
\lim_{t\to 0}\myint{B_\gr(z)}{}u(t,x) dx=+\infty\quad{\rm for\ any} \  \, \gr>0.
\ee
 They also highlighted the existence of a critical exponent $p_c=1+\frac 2N$, which plays a crucial role in the fine analysis of the initial trace. For example they obtained that if $p$ is {\it subcritical}, i.e. $1<p<p_c$,  $(\ref{In-4})$ can be sharpened in the form
\bel{In-8}
c_2(p,N)\leq \liminf_{t\to 0}\,t^{\frac{1}{p-1}}u(z,t)\leq \limsup_{t\to 0}\,t^{\frac{1}{p-1}}u(z,t)\leq c_1(p),
\ee
for some positive constants $c_1(p)>c_2(p,N)$. Furthermore they proved that for any couple $(\CS,\gn)$, where $\CS$ is a closed subset of $\BBR^N$ and $\gn$ a nonnegative Radon measure on $\CR=\BBR^N\setminus\CS$, there exists a unique nonnegative solution $u$ of $(\ref{In-6})$ with the initial trace $Tr(u)=(\CS,\gn)$. The {\it supercritical} case $p\geq p_c$ turned out to be much more delicate and was finally elucidated in a series of works by Marcus and V\'eron \cite{MV5} and Gkikas and V\'eron \cite{GV} following some deep ideas introduced by Marcus and V\'eron in \cite{MV4} and Marcus \cite{M} for solving similar questions dealing with semilinear elliptic equations. Al Sayed and V\'eron in \cite {AV} extended  the subcritical analysis performed in \cite{MV1} to the non-autonomous equation
\bel{In-9}
\prt_tu-\Gd u+t^\gb u^p=0\quad{\rm in }\; \;Q_\infty,
\ee
with $\gb>-1$ and $p>1$. Note that the choice $\gb>-1$ is natural otherwise the initial trace would be essentialy zero as it can be verified with the equation without absorption. \medskip

The main difficulty to extend some of the previous results dealing with $(\ref{In-6})$ and $(\ref{In-9})$ comes from the fact that the fractional Laplacian is a non-local operator. A more precise characterization of the conditional singular set needs additional assumptions on $u$ or on $f$. We define the {\it singular set }$\CS_u$ of $u$ by
\bel{In-9sup}\displaystyle
 \CS_u=\left\{z\in\BBR^N:\limsup_{t\to 0}\myint{B_\gr(z)}{}u(t,x) dx=+\infty\ \, {\rm for\ any}\ \,  \gr>0\right\}.
\ee
This set is closed and it follows from Theorem A that $\CS_u\subset \tilde \CS_u$.  The {\it initial trace} is the couple $Tr(u):=(\CS_u,\gn)$. This initial trace can also be seen as an outer regular Borel measure with regular part (or Radon part) $\gn$ and singular part $\CS_u$.
When $s=1$ then $Tr(u)=Tr_c(u)$ because the set $ \CS_u$ is also characterized as the set of $z\in \BBR^N$ where 
 \bel{In-9sup1}\displaystyle
\myint{0}{\frac T2}\myint{B_\gr(z)}{}f(t,x,u)dxdt=\infty\ \, {\rm for\ any}\ \,  \gr>0.
\ee
When $0<s<1$ and no extra assumption on $f$ are made, $Tr(u)$ could be different from $Tr_c(u)$. 
\medskip

\nind{\bf Theorem B} {\it Assume that $u$ is a nonnegative solution of (\ref{In-1}). If $u\in L^1(0,T;\CL^s(\BBR^N))$, then
$\CS_u=\tilde\CS_u$ and more precisely for any $z\in\CS_u$,
\bel{In-13}
\lim_{t\to 0}\myint{B_\gr(z)}{}u(t,x)\, dx=+\infty\ \  {\rm for\ any}\ \,  \gr>0.
\ee
}\smallskip

The above assumption on $u$ can be verified when the absorption is strong and the singular set is compact. Another type of characterization of the singular set  needs the following assumptions on $f$:  $f(t,x,u)$ satisfies $f(t,x,0)=0$ and
\bel{In-10}
0\leq f(t,x,u)\leq t^\gb g(u)\quad\forall\, (t,x,u)\in \BBR_+\ti\BBR^N\ti\BBR_+,
\ee
where $\R_+=[0,+\infty)$,   $\gb>-1$, $g$ is  nondecreasing, continuous and verifies the subcritical growth assumption,
\bel{In-11}
\myint{1}{\infty}g(s)s^{-1-p_\gb^*}ds<+\infty,
\ee
with
\bel{In-12}
p_\gb^*=1+\myfrac{2s(1+\gb)}{N}.
\ee
The role of the subcritical growth assumption (\ref{In-1}) has been highlighted in \cite{CVW} as the natural condition to solve 
the initial value problem with a bounded positive Radon measure for equation (\ref{In-1}) (see Section 2.2).  \medskip

\nind{\bf Theorem C} {\it Assume (\ref{In-10}) and either (\ref{In-11}) holds if $-1<\gb\leq 0$, or
\bel{In-11^*}
\myint{1}{\infty}g(s)s^{-2-\frac{2s}{N}}ds<+\infty,
\ee
 if $\gb>0$. If $u$ is a nonnegative solution of (\ref{In-1}) with initial trace $(\CS_u,\gn_u)$. If $\CS_u\neq\emptyset$ and $z\in\CS_u$, then (\ref{In-13}) holds. More precisely $u\geq  u_{z,\infty}$ where $ \displaystyle u_{z,\infty}=\lim_{k\to\infty} u_{k\gd_z}$ and $ u_{k\gd_z}$ is the solution of
\bel{In-14}\BA{lll}
\partial_t  u + (-\Delta)^s  u+t^\gb g( u)=0\quad  {\rm in}\quad Q_\infty\\[1.2mm]
\phantom{----..---\, }
 u(0,.)=k\gd_z.
\EA\ee
}\medskip

The existence and uniqueness of  solutions to (\ref{In-14}) follow from \cite[Th 1.1]{CVW}. If $g:\BBR\rightarrow \BBR_+$ is nondecreasing and satisfies that
\bel{In-14-1}\BA{lll}
 \mathcal{G}(t):=\myint{t}{\infty}\myfrac{ds}{g(s)}<+\infty\quad {\rm for }\ \, t>0,
\EA\ee
and if $\gb>-1$, denote $U(t)=\mathcal{G}^{-1}(\myfrac{t^{\gb+1}}{\gb+1})$,
where $\mathcal{G}^{-1}$ is the inverse function of $\mathcal{G}$, then  the function $U$ verifies that
\bel{In-14-2}\BA{lll}
\myint{U(t)}{\infty}\myfrac{ds}{g(s)}=\myfrac{t^{\gb+1}}{\gb+1},
\EA\ee
and defines as the maximal solution of the ODE
\bel{In-14-3}\BA{lll}
\prt_tU+t^\gb g(U)=0\quad{\rm on }\;\,\BBR^*_+\quad{\rm satisfying }\;\ U(0)=+\infty.
\EA\ee

\medskip

\nind{\bf Theorem D} {\it Assume that $f(t,x,r)\geq t^\gb g(r)$, where $\gb>-1$ and $g$ satisfies (\ref{In-14-1}). If  $u$ is a nonnegative solution of (\ref{In-1}) belonging to $L^1_{loc}(0,T;\CL^s(\BBR^N))$, then
\bel{In-13sup}
u(t,x)\leq U(t),\quad\forall\, (t,x)\in Q_\infty.
\ee
Furthermore, if $g$ satisfies
\bel{In-14-4}
\myint{1}{\infty}\myfrac{sds}{g(s)\left(\myint{s}{\infty}\myfrac{d\gt}{g(\gt)}\right)^{\frac{\gb}{\gb+1}}}<+\infty,
\ee
then $\CS_u=\tilde\CS_u$ and (\ref{In-13}) holds for any $z\in \CS_u$.

}\medskip

\medskip

\nind{\bf Theorem E} {\it Assume that $f(t,x,s)= t^\gb g(s)$, where $\gb>-1$ and $g$ satisfies (\ref{In-14-1}), is nondecreasing and is locally Lipschitz continuous. If  $u$ is a nonnegative solution of (\ref{In-1}) belonging to $L^1_{loc}(0,T;\CL^s(\BBR^N))$, then either
\bel{In-14-5}
u(t,x)< U(t),\quad\forall\, (t,x)\in Q_\infty,
\ee
or
\bel{In-14-6}
u(t,x)= U(t),\quad\forall\, (t,x)\in Q_\infty.
\ee
}\medskip


In the second part of this paper we study in detail the initial trace problem for the equation
\bel{In-15}
\prt_tu+(-\Gd)^s u+t^\gb u^p=0\quad{\rm in }\quad Q_\infty,
\ee
 when $s\in(0,1)$, $\beta>-1$ and $p\in(1,p_\gb^*)$. A second critical value of $p$ appears
 \bel{In-16}
p_\gb^{**}=1+\myfrac{2s(1+\gb)}{N+2s}.
\ee
Actually, if $u_k:=u_{k\gd_0}$ is unique solution to
\begin{equation}\label{In-16'}
\BA{lll}
\partial_t u + (-\Delta)^s  u+ t^\beta u^p=0\quad & {\rm in}\quad Q_\infty,\\[1 mm]
\phantom{ \partial_t u + -\Delta^s  t^\beta u^p }
u(0,\cdot)=k\delta_{0}\quad & {\rm in}\quad \R^N,
\EA
\end{equation}
it is proved in \cite{CVW} that $\displaystyle u_\infty=\lim_{k\to\infty} u_k$ is very different according $1<p< p_\gb^{**}$ or $p_\gb^{**}<p<p_\gb^{*}$. Precisely,
\nind (i) if $1<p< p_\gb^{**}$, then
 \begin{equation}\label{In-17}
   u_\infty(t,x)=U_{p,\gb}(t):=\left(\frac{1+\beta}{p-1}\right)^{\frac1{p-1}} t^{-\frac{1+\beta}{p-1}}.
 \end{equation}
 {\it The absorption is dominant}, as if $s=0$. \smallskip

\nind (ii) If $p_\gb^{**}<p< p_\gb^{*}$, then
 \begin{equation}\label{In-18}
   u_\infty(t,x)=V(t,x):=t^{-\frac{1+\beta}{p-1}}v\left(\frac{x}{t^{\frac{1}{2s}}}\right),
 \end{equation}
 where $v$ is the minimal positive solution of
  \begin{equation}\label{In-19}\BA {lll}
  (-\Gd)^s v-\myfrac{1}{2s}\nabla v\cdot\eta-\myfrac{1+\gb}{p-1}v+v^p=0\quad{\rm in }\quad \R^N,\\[2.2mm]
  \phantom{------,\,\,\cdot}
  \displaystyle
\lim_{|\eta|\to\infty}|\eta|^{\frac{2s(1+\gb)}{p-1}}v(\eta)=0.
 \EA\end{equation}
 The function $V$  is called the very singular solution of (\ref{In-15}). In this case {\it the diffusion is dominant}, as when $s=1$. 
  \smallskip
  
  We first prove the following result which complements Theorem C in the case where $\gb>0$. The proof is delicate and uses a 
  form of parabolic Harnack inequality valid for solutions of (\ref{In-15}).\medskip
  
  \nind{\bf Theorem F} {\it Assume $\gb>-1$, $1<p<p^*_\gb$ and $u$ is a nonnegative solution of (\ref{In-15}) with initial trace $(\CS_u,\gn_u)$. If $\CS_u\neq\emptyset$ and $z\in\CS_u$ then $u\geq u_{z,\infty}$.}\medskip

We observe that
 $\tilde\mathcal{S}_{u_\infty}=\mathcal{S}_{u_\infty}=\{0\}$ when $p_\gb^{**}<p\leq p_\gb^{*}$
 and   $\tilde\mathcal{S}_{u_\infty}=\mathcal{S}_{u_\infty}=\R^N$ when $1<p<p_\gb^{**}$. Notice that the case $p=p_\gb^{**}$ remained unsolved in \cite{CVW}.
In this paper, we prove that $\mathcal{S}_{u_\infty}=\R^N$ also for $p=p_\gb^{**}$. Our main result concerning (\ref{In-15}) is the following.\medskip

\nind {\bf Theorem G  }{\it Let $u$ be a positive solution of (\ref{In-15}).\smallskip

\nind (i) If  $ p\in(1,p_\gb^{**}]$  and  $\mathcal{S}_{u}\not=\emptyset$. Then $\mathcal{S}_u=\R^N$ and $u\geq U_{p,\gb}$. If we assume moreover that
$u\in L^1_{loc}(0,\infty;\CL^s(\BBR^N))$, then $u=U_{p,\gb}$.\smallskip

\nind (ii) If  there exists $\gk\in [1,N]\cap\BBN$ such that $p\in(1,p_\gb^{*})\cap\left(1,1+\frac{2s(1+\beta)}{\gk+2s}\right]$
and $\CS_u$ contains an affine plane $\mathcal{L}$ of codimension $\gk$.
 Then the conclusions of (i) hold}. \medskip

If $\gk=N$, (ii) is just (i). Note that if $0<s\leq \frac{1}{2}$ or if $\gk\geq N-2$, then $(p_\gb^{**},p_\gb^{*})\cap\left(p_\gb^{**},1+\frac{2s(1+\beta)}{\gk+2s}\right]=(p_\gb^{**},p_\gb^{*})$, while, if $\frac{1}{2}< s<1$ and
$\gk=N-1$, then $(p_\gb^{**},p_\gb^{*})\cap\left(p_\gb^{**},1+\frac{2s(1+\beta)}{\gk+2s}\right]=\left(p_\gb^{**},1+\frac{2s(1+\beta)}{N-1+2s}\right]$. \smallskip

Conversely, given  a closed set of $\CS\subset\R^N$ and a nonnegative Radon measure on $\nu$ on  $\CR=\R^N\setminus \mathcal{S}$, we study the existence of solution of (\ref{In-15}) with a given initial trace
$Tr_c(u)=Tr(u)=(\mathcal{S},\nu)$,  that is a solution of the following problem
\begin{equation}\label{In-22}
 \arraycolsep=1pt
\begin{array}{lll}
\partial_t u + (-\Delta)^s  u+ t^\beta u^p=0\quad & {\rm in}\quad Q_\infty\\[1mm]
\phantom{-------}
Tr(u)=(\mathcal{S},\nu).
\end{array}
\end{equation}
This means that $u$ is a classical solution of the equation in $Q_\infty$ and that (\ref{In-4}) and (\ref{In-13sup}) hold. By Theorem G  any closed set cannot be the singular part of the initial trace of a positive solution of (\ref{In-15}) if $p$ is too small (diffusion effect) or if $p$ is too large. In the same sense any positive bounded Radon measure $\gn$ cannot be the regular part of the initial trace of a positive solution of (\ref{In-15}) since condition (\ref{In-11}) is equivalent to $p<p^*_\gb$. However this condition is restrictive and there exist several sufficient conditions linking $\gn$, $s$, $\gb$ and $p$. Hence we say that a nonnegative bounded measure $\gn$ is an {\it admissible measure} if the initial value problem
\begin{equation}\label{In-23}
 \arraycolsep=1pt
\begin{array}{lll}
\partial_t u + (-\Delta)^s  u+ t^\beta u^p=0\quad & {\rm in}\quad Q_\infty\\[1mm]
\phantom{------,,\,\,}
u(0,.)=\nu,
\end{array}
\end{equation}
admits a solution $u_\gn$, always unique, and it is a {\it good measure} if it is stable in the sense that if $\gn$ is replaced by $\gn\ast\gr_n$ for some sequence of mollifiers, then $u_{\gn\ast\gr_n}$ and $t^\gb u^p_{\gn\ast\gr_n}$ converges to $u_{\gn}$ and $t^\gb u^p_{\gn}$ respectively in $L^1(Q_T)$. We denote by $H_s$ is the  kernel in $\BBR^*_+\times\R^N$ associated to $(-\Gd)^s$. It is expressed by
   \begin{equation}\label{In-24}\BA {lll}
H_s(t,x)=\myfrac{1}{t^{\frac{N}{2s}}}\tilde H_s\left(\frac{x}{t^{\frac{1}{2s}}}\right)\quad{\rm where }\quad\tilde H_s(x)=\myint{\BBR^N}{}e^{ix.\xi-\abs{\xi}^{2s}}d\xi.
 \EA\end{equation}
 It  is proved in \cite{CKS}, \cite {BSV} that $H_s$ satisfies the following two-side estimate,
 \bel{cks}
\myfrac{c^{-1}_{3}t}{t^{1+\frac{N}{2s}}+|x|^{N+2s}}\leq H_s(t,x)\leq\myfrac{c_{3}t}{t^{1+\frac{N}{2s}}+|x|^{N+2s}}\quad\forall\, (t,x)\in\BBR_+^*\ti\BBR^N.
 \ee
The associated potential $\mathbb{H}_s[\nu]$ of a bounded Radon measure $\gn$ in $\BBR^N$  is defined by
$$\mathbb{H}_s[\nu](t,x)=\int_{\R^N} H_s(t,x-y)d\nu(y).$$
We first prove that a nonnegative bounded measure with Lebesgue decomposition $\gn=\gn_0+\gn_s$, where $\gn_0\in L^1(\BBR^N)$ and $\gn_s$ is singular with respect to the $N$-dim Lebesgue measure is a good measures if $t^\gb(\mathbb{H}_s[\nu])^p\in L^1(Q_1)$. The expression of admissibility for a Radon measure needs the introduction of Bessel capacities which are presented in Section 4.1. Our main existence result of solutions for (\ref{In-23}) is the following.\medskip

\nind{\bf Theorem H } {\it Let $N\geq 1$, $p>1$ and $-1<\gb<p-1$. Then a nonnegative bounded measure $\gn$ in $\BBR^N$ is an admissible measure if and only if $\gn$ vanishes on Borel subsets of $\BBR^N$ with zero  ${\rm cap}^{\BBR^N}_{\frac{2s(1+\gb)}{p},p'}$-Bessel  capacity.}

\medskip
Concerning problem (\ref{In-22}) we have the following general result.\medskip

\nind {\bf Theorem I  }{\it Assume that $N\geq 1$ and $p>1+\frac{2s (1+\gb)}{1+2s}$. If $\CS$ is a closed subset of $\BBR^N$ such that $\CS=\overline{{\rm int}\,\CS}$ and $\gn$ is a nonnegative Radon measure on $\CR=\CS^c$ such that for any compact set $K\subset\CR$,
$\chi_{_K}\gn$ is an admissible measure. Then problem (\ref{In-22}) admits a solution.
}
\medskip

It is interesting to compare this result with \cite[Th. 4.11]{MV1} dealing with the case  $s=1$.  It is proved there that for any closed set satisfying a non-thinness condition (expressed later on in terms of Bessel capacity \cite{GV}) but always fulfilled when $1<p<1+\frac 2N$ and any nonnegative admissible Radon measure, problem (\ref{In-22}) admits a solution. There the condition $p>1+\frac{2s (1+\gb)}{1+2s}$ has no counterpart when $s=1$. Theorem L shows that this condition is fundamental in order to have existence without condition at infinity, even in the case  where $\CS=\emptyset$ and  $\gn$ is a mere $L^1_{loc}(\BBR^N)$ function.
In some particular cases,  the existence of a solution to (\ref{In-22}) with no extra condition on $\CS$ or $\gn$ can be proved as the next results show it.\medskip

\nind {\bf Theorem J }{\it Assume that $\gb>-1$, $p>1$ and  one of the following assumptions is fulfilled:\smallskip

\nind (i) either $N=1$ and $1+\frac{2s(1+\gb)}{1+2s}<p<1+2s(1+\gb)$,\smallskip

\nind (ii) or\quad  $N=2$, $\frac{1}{2}<s<1$ and $1+\frac{2s(1+\gb)}{1+2s}<p<1+s(1+\gb)$.\smallskip

\nind
Then for any closed set $\CS\subset\R^N$ and any
nonnegative measure $\gn$ in $\CR=\CS^c$, there exists a nonnegative solution $u$ to (\ref{In-22}). }\medskip

As a consequence of the previous results we obtain existence with initial data measure in $\BBR^N$ of solutions without condition at infinity in the spirit of Brezis classical result \cite{Brez}.\medskip

\nind{\bf Corollary K} {\it Let $N\geq 1$, $\gb>-1$ and $p>1$. \smallskip

\nind (a) If $p>1+\frac{2s (1+\gb)}{1+2s}$ and $-1<\gb<p-1$, then there exists a positive solution $u\in L^1_{loc}(0,\infty;\CL^s(\BBR^N))$ to problem (\ref{In-23}) for any nonnegative Radon measure $\gn$ in $\BBR^N$ if and only if for any $n\in\BBN^*$, $\chi_{_{B_n}}\gn$ is an admissible measure.\smallskip

\nind (b) If one of the assumptions (i) or (ii) of Theorem J is fulfilled, then for any nonnegative Radon measure $\gn$ in $\BBR^N$ there exists a positive solution to problem (\ref{In-23}). \medskip
}

Conditions (i) or (ii) of Theorem J are essentially necessary for unconditional existence, since we have the following result.
 \medskip

\nind {\bf Theorem L }{\it Assume that $\gb>-1$, $p>1$ and $1<p<1+\frac{2s(1+\gb)}{1+2s}$. If $\gf\in L^1_{loc}(\BBR^N)$ is a nonnegative function which satisfies
   \begin{equation}\label{In-25}\BA {lll}\displaystyle
\lim_{\abs x\to\infty}\myfrac{\gf(x)}{|x|^{\frac{2s(\gb+1)}{\gb+2-p}}}=+\infty,
 \EA\end{equation}
 then the sequence of solutions $\{u_n\}$ of (\ref{In-23}) with initial data
 $\gn=\gn_n=\inf\{\gf,\gf_n\}$, where $\gf_n=\inf\{\gf(z):|z|\geq n\}$ is increasing and converges to
 $U(t)$.} \medskip

This implies that there exists no solution of (\ref{In-23}) with initial data $\gf$. Notice that $\gb+2-p$ is positive for $p<1+\frac{2s(1+\gb)}{1+2s}$. For the mere heat equation a theory of maximal growth for admissibility growth of initial data has been developed in \cite {Wi} and for the fractional heat equation in \cite {BSV}. In both cases the representation formulas play an important role. 
For equations with potential a phenomenon of instantaneous blow-up is proved \cite{BG} for solution of 
   \begin{equation}\label{In-26}\BA {lll}\displaystyle
\prt_tu-\Gd u-V(x)u=0\quad  {\rm in}\quad Q_\infty,
 \EA\end{equation}
when $V\sim c\abs x^{-2}$, for any nonnegative initial data. This phenomenon of instantaneous blow-up has been recently highlighted in \cite{ShVe} for the the semilinear equation
   \begin{equation}\label{In-27}\BA {lll}\displaystyle
\prt_tu-\Gd u+u\left(\ln(u+1)\right)^\ga=0\quad  {\rm in}\quad Q_\infty,
 \EA\end{equation}
 when $1<\ga<2$. It is shown there  that the limit of a sequence of solutions with fast growing initial data is the maximal solution of $u'+u\left(\ln(u+1)\right)^\ga=0$ on $(0,\infty)$.
\bigskip\bigskip

 \noindent{\bf Acknowledgements} \smallskip
 
 \nind H. Chen  is supported by NNSF of China, No: 11726614, 11661045,   by the Jiangxi Provincial Natural Science Foundation, No: 20161ACB20007. The authors are grateful to the referee for the careful reading of the manuscript which enabled several  improvements in the presentation of this work.


\setcounter{equation}{0}
\section{Initial trace with general nonlinearity}

\subsection{Existence of an initial trace}
\nind{\bf Proof of Theorem A}. For any bounded domain $\gw\subset\R^N$, we denote by $C^2_{0}(\overline\gw)$ the space of functions $\xi:\R^N\rightarrow\R$ which are $C^2$ and have compact support in $\overline\gw$. We always assume that  $N\geq1$ and $0<s<1$. Let $\gf_\gw$ be the first eigenfunction of $(-\Gd)^s$ in $H^{s}_0(\gw)$, with corresponding eigenvalue $\gl_\gw>0$, i.e. the solution of
\bel{II-1}\BA {lll}
(-\Gd)^s\gf_\gw=\gl_\gw\gf_\gw&\quad{\rm in}\quad \gw\\
\phantom{(-\Gd)^s}
\gf_\gw=0&\quad{\rm in}\quad \gw^c.
\EA\ee
Existence and basic properties of the eigenfunctions can be found in \cite{BK}, \cite{BBKRSZ}. We normalize $\gf_\gw$ by $\sup \gf_\gw=1$. We say that $\gw$ is of class {\bf E. S. C.} if it satisfies the {\it exterior sphere condition}. It is known by \cite[Prop 1.1]{RO-S} that $\gf_\gw(x)\leq c(\dist (x,\prt\gw))^s$ in $\gw$, and there exists $q>2$ such that $\gf^q_\gw\in C^2_{0}(\overline\gw)$.  We denote by $K_\gr(z)$ the open cube with sides parallel to the axis of center $z\in\BBR^N$ and length sides $\gr>0$, and $K_1:=K_1(0)$.
Then
$$\gf_{K_\gr(z)}(x)=\gf_{K_1}\left(\frac{x-z}{\gr}\right)\quad{\rm and}\quad \gl_{K_\gr(z)}=\frac{\gl_{K_1}}{\gr^{2s}}.$$
The next lemma is an improvement of \cite[Lemma 2.3]{CV2}.
\blemma{prod} Let $q\in\BBN\cap [2,\infty)$ and $\gz\in C^2_{0}(\overline\gw)$, $\gz\geq 0$, then
\bel{II-2}\BA {lll}
(-\Gd)^s\gz^q(x)=q\gz^{q-1}(x)(-\Gd)^s\gz(x)-\myfrac{a_{N,s}}{q(q-1)}\myint{\BBR^N}{}\myfrac{\gz^{q}(y)-\gz^{q}(x)-q(\gz(y)-\gz(x))\gz^{q-1}(y)}{\abs{x-y}^{N+2s}}dy\\[4mm]
\phantom{(-\Gd)^s\gz^q(x)}
\geq q\gz^{q-1}(x)(-\Gd)^s\gz(x)-\myfrac{a_{N,s}}{q}\gz^{q-2}(x)\myint{\BBR^N}{}\myfrac{\left(\gz(y)-\gz(x)\right)^2}{\abs{x-y}^{N+2s}}dy.
\EA\ee
\es
\Proof From  \cite[Lemma 2.3]{CV2}, we know that
$$\BA {lll}(-\Gd)^s\gz^q(x)=q\gz^{q-1}(x)(-\Gd)^s\gz(x)-q(q-1)a_{N,s}\myint{\BBR^N}{}\left(\myint{\gz(x)}{\gz(y)}(\gz(y)-t)t^{q-2}dt\right)\myfrac{dy}{\abs{x-y}^{N+2s}}.
\EA$$
By integration by parts, we obtain that
$$\BA {ll}\myint{\gz(x)}{\gz(y)}(\gz(y)-t)t^{q-2}dt=\myfrac{1}{q(q-1)}\left(\gz^q(y)-\gz^q(x)-q(\gz(y)-\gz(x))\gz^{q-1}(x)\right)\\[2mm]
\phantom{\myint{\gz(x)}{\gz(y)}----}
=\myfrac{\gz(y)-\gz(x)}{q(q-1)}\left[\gz^{q-1}(y)+\gz^{q-2}(y)\gz(x)+...+\gz(y)\gz^{q-2}(x)-(q-1)\gz^{q-1}(x)\right].
\EA$$
Since for any $a,b\geq 0$
$$\BA{ll}
b^{q-1}+ab^{q-2}+...a^{q-2}b-(q-1)a^{q-1}\\[1mm]\phantom{----}
=b^{q-1}-a^{q-1}+a(b^{q-2}-a^{q-2})+a^2(b^{q-3}-a^{q-3})+...+a^{q-2}(b-a)\\[1mm]\phantom{----}
=(b-a)\left[\left(b^{q-2}+ab^{q-3}+...+a^{q-2}\right)+a\left(b^{q-3}+ab^{q-4}+...+a^{q-3}\right)+...+a^{q-2}\right]
\\[1mm]\phantom{----}
\geq (q-1)(b-a)a^{q-2},
\EA$$
we obtain (\ref{II-2}).\qeda\medskip

\nind\Remark By the mean value theorem, we see that there exists $m_\gz\in\{z=\gz(w):w\in\BBR^N\}$ such that
\bel{II-3}\BA {lll}
L(\gz^q):=\myfrac{a_{N,s}}{q(q-1)}\myint{\BBR^N}{}\myfrac{\gz^q(y)-\gz^q(x)-q(\gz(y)-\gz(x))\gz^{q-1}(x)}{\abs{x-y}^{N+2s}}dy\\[4mm]
\phantom{L(\gz^q):}
=\myfrac{a_{N,s}}{2}m_\gz^{q-2}\myint{\BBR^N}{}\myfrac{\left(\gz(y)-\gz(x)\right)^2}{\abs{x-y}^{N+2s}}dy.
\EA\ee

\bprop{loc1}  Assume that $f$ satisfies (\ref{In-2}) and $u$ is a nonnegative solution of (\ref{In-1}) such that $u(t,\cdot)\in \CL^s(\BBR^N)$ for all $t\in (0,T)$. If $f(\cdot,\cdot,u)\in L^1(Q^T_\gw)$ for some  bounded domain $\gw\subset\R^N$ of class {\bf E. S. C.} and $T>0$. Then there exists $\ell_{\gw}\geq 0$ such that
\bel{II-4}\BA {lll}\displaystyle
\lim_{t\to 0}\myint{\gw}{}u(t,x)\gf^q_\gw(x)dx=\ell_{\gw}.
\EA\ee
Furthermore, we have that
\bel{II-5}\BA {lll}\displaystyle
\ell_{\gw}+\myfrac{a_{N,s}}{q}\myint{0}{T}\myint{\BBR^N}{}u(s,x)\gf^{q-2}_\gw(x)\left(\myint{\BBR^N}{}\myfrac{\left(\gf_\gw(y)-\gf_\gw(x)\right)^2}{\abs{x-y}^{N+2s}}dy\right)dx\\[4mm]
\phantom{----------}
\leq e^{q\gl_\gw T}X(T)+\myint{0}{T}\myint{\gw}{}f(s,x,u)\gf^q_\gw(x)e^{q\gl_\gw s}dxdt.
\EA\ee
\es
\Proof Since $\gf^q_\gw\in C^2_{0}(\overline\gw)$, there holds
\bel{II-6-}\myfrac{d}{dt}\myint{\gw}{}u(t,x)\gf^q_\gw(x)dx+\myint{\BBR^N}{}u(t,x)(-\Gd)^s\gf^q_\gw(x)dx+\myint{\gw}{}f(t,x,u)\gf^q_\gw(x)dx=0.
\ee
Set
$$X(t)=\myint{\gw}{}u(t,x)\gf^q_\gw(x)dx,
$$
then
\bel{II-6}\BA {lll}\myfrac{d}{dt}\left(e^{q\gl_\gw t}X(t)-\myint{t}{T}\myint{\gw}{}f(s,x,u)\gf^q_\gw(x)e^{q\gl_\gw s}dxds\right)=
e^{q\gl_\gw t}\myint{\BBR^N}{}L(\gf^q_\gw)(x)dx\geq 0.
\EA\ee
This implies that $\lim_{t\to 0}X(t)=\ell_\gw$ exists and
$$\BA {lll}
\ell_\gw+\myint{0}{T}\myint{\BBR^N}{}L(\gf^q_\gw)(x)e^{q\gl_\gw s}u(s,x)dxds+\myint{0}{T}\myint{\gw}{}f(s,x,u)\gf^q_\gw(x)e^{q\gl_\gw s}dxds= e^{q\gl_\gw T}X(T),
\EA$$
which implies (\ref{II-5}) by \rlemma{prod}.\qeda\medskip

As an immediate consequence we have,

\bcor{loc2}  Under the assumptions of Theorem A, $u\in L^1(Q_G^T)$ for any compact set $G\subset\CR_u$.
\es
\medskip

The proof of Theorem A is completed by the following statement:
\bprop{tr} There exists a nonnegative Radon measure $\mu_u$ on $\CR_u$ such that for any $\gz\in C^2_0(\CR_u)$, there holds
\bel{II-7}\BA {lll}\displaystyle
\lim_{t\to 0}\myint{\CR_u}{}u(t,x)\gz(x)dx=\myint{\CR_u}{}\gz d\gm_u.
\EA
\ee
\es
\Proof Let $\gz\in C^2_0(\CR_u)$ with support $K$ and let $G$ be an open subset containing $K$ such that $\prt G$ is smooth and $\overline G$ is a compact subset of $\CR_u$ and assume $0\leq\gz\leq 1$. We put
$$Y(t)=\myint{\BBR^N}{}u(t,x)\gz(x)dx=\myint{G}{}u(t,x)\gz(x)dx,
$$
and
$$\myint{\BBR^N}{}f(t,x,u)\gz(x) dx=\myint{G}{}f(t,x,u)\gz(x) dx,
$$
then
$$Y'(t) +\myint{\BBR^N}{}u(t,x)(-\Gd)^s\gz(x) dx+\myint{G}{}f(t,x,u)\gz(x) dx=0.
$$
Since $\gz\geq 0$, we have that
\bel{II-7'}\BA {lll}\myfrac{1}{a_{N,s}}\myint{\BBR^N}{}\!u(t,x)(-\Gd)^s\gz(x) dx\\[4mm]
\phantom{-\myfrac{1}{a_{N}}\myint{\BBR^N}{}\!u(t,x) dx}
=\!\myint{\BBR^N}{}\!u(t,x)\myint{G}{}\myfrac{\gz(x)-\gz(y)}{\abs{x-y}^{N+2s}}dy dx
+\!\!\myint{G}{}\!u(t,x)\gz(x)\myint{G^c}{}\myfrac{dy}{\abs{x-y}^{N+2s}} dx\\[4mm]
\phantom{-\myfrac{1}{a_{N}}\myint{\BBR^N}{}\!u(t,x) dx}
=\!\myint{G}{}\!u(t,x)\myint{G}{}\myfrac{\gz(x)-\gz(y)}{\abs{x-y}^{N+2s}}dy dx-\!\myint{G^c}{}\!u(t,x)\myint{G}{}\myfrac{\gz(y)}{\abs{x-y}^{N+2s}}dy dx\\[4mm]
\phantom{\myint{\BBR^N}{}\!u(t,x)(-\Gd)^s\gz(x) dx------------}
+\!\!\myint{G}{}\!u(t,x)\gz(x)\myint{G^c}{}\myfrac{dy}{\abs{x-y}^{N+2s}} dx\\[4mm]
\phantom{-\myfrac{1}{a_{N}}\myint{\BBR^N}{}\!u(t,x) dx}
\leq \!\myint{G}{}\!u(t,x)\myint{G}{}\myfrac{\gz(x)-\gz(y)}{\abs{x-y}^{N+2s}}dy dx+\!\!\myint{G}{}\!u(t,x)\gz(x)\myint{G^c}{}\myfrac{dy}{\abs{x-y}^{N+2s}} dx
\EA\ee
If we define the regional fractional Laplacian of order $s$ relative to $G$ by
$$(-\Gd)^s_G\gz(x):=a_{N,s}\myint{G}{}\myfrac{\gz(x)-\gz(y)}{\abs{x-y}^{N+2s}}dy,
$$
then the right-hand side of $(\ref{II-7'})$ is bounded from above by
$$\Gl(t)=\displaystyle\left(\norm{(-\Gd)^s_G\gz}_{L^\infty}+\max_{x\in K}\myint{G^c}{}\myfrac{dy}{\abs{x-y}^{N+2s}}\right)
\myint{G}{}\!u(t,x) dx,
$$
since $\gz$ is $C^2$ with support in $K\subset G\subset\overline G\Subset\CR_u$. By \rcor{loc2}, $\Gl\in L^1(0,T)$. Because
\bel{II-8}
\myfrac{d}{dt}\left(Y(t) -\myint{t}{T}\left(\Gl(s)+\myint{G}{}f(t,x,u)\gz(x) dx\right)ds\right)\geq 0,
\ee
and
\bel{II-9}\myint{0}{T}\left(\Gl(s)+\myint{G}{}f(t,x,u)\gz(x) dx\right)ds<\infty
\ee
Combining (\ref{II-8}) and (\ref{II-9}) we infer that the following limit exists
\bel{II-10}\displaystyle
\lim_{t\to 0}Y(t)=\lim_{t\to 0}\myint{G}{}u(t,x)\gz(x)dx:=\tilde\gm_u(\gz).
\ee
By replacing $\gz$ by $\norm\gz_{L^\infty}\gz$ we can drop the condition  $\gz\leq 1$. with support in $K$, then
\bel{II-11}\displaystyle
0\leq \lim_{t\to 0}\myint{G}{}u(t,x)\gz(x)dx=\tilde\gm_u(\gz)\leq \ell_G\norm\gz_{L^\infty}.
\ee

Next we assume that $\gz\in C_0(\CR_u)$ is nonnegative, with support $K\subset G\subset\overline G\Subset\CR_u$, then there exists an increasing sequences $\{\gz_n\}\subset C^2_0(\CR_u)$ of nonnegative functions smaller than $\gz$ which converges to $\gz$ uniformly (take for example $\gz_n=(\gz-n^{-1})_+\ast\gr_n$ for some sequence of mollifiers $\{\gr_n\}$ with supp$(\gr_n)\subset B_{n^{-2}}$). The sequence $\{\tilde\gm_u(\gz_n)\}$ is increasing and bounded from above by $M\ell_G\sup_G\gz$. Hence it is convergent and its limit, still denoted by $\tilde\gm_u(\gz)$ is independent of the sequence $\{\gz_n\}$. We can also consider a uniform approximation of $\gz$ from above in considering $\gz'_n=(\gs_n+\gz)\ast\gr_n$,  where $\gs_n=n^{-1}\chi_{_{K_n}}$ and $K_n=\{x\in\BBR^N:\dist (x,K)\leq n^{-1}\}$. Actually,
\bel{II-12}\displaystyle
\tilde\gm_u(\gz)=\sup\{\tilde\gm_u(\eta):\eta\in C^2_0(\CR_u),0\leq\eta\leq\gz\}=\inf\{\tilde\gm_u(\eta'):\eta'\in C^2_0(\CR_u),\gz\leq\eta'\}.
\ee
This implies that for all $\eta$ and $\eta'$ belonging to $C^2_0(\CR_u)$ such that $\eta\leq\gz\leq\eta'$, we have that
\bel{II-13}\displaystyle
\tilde\gm_u(\eta)\leq\liminf_{t\to 0}\myint{\CR_u}{}u(t,x)\gz(x)dx\leq \limsup_{t\to 0}\myint{\CR_u}{}u(t,x)\gz(x)dx\leq \tilde\gm_u(\eta').
\ee
Combined with (\ref{II-12}) we infer the existence of the limit and
\bel{II-14}\displaystyle
\lim_{t\to 0}\myint{\CR_u}{}u(t,x)\gz(x)dx=\tilde\gm_u(\gz).
\ee
Finally, if $\gz\in C_0(\CR_u)$ is a signed function we write $\gz=\gz_+-\gz_-$ and $\gm_u(\gz)=\tilde\gm_u(\gz_+)-\tilde\gm_u(\gz_-)$. Hence $\gm_u$ is a positive Radon measure on $\CR_u$, and (\ref{II-7}) follows from (\ref{II-14}) with $\gz$ replaced by
$\gz_+$ and $\gz_-$.\qeda
\medskip

\blemma{test}
Assume that $G\subset\BBR^N$ is a bounded smooth domain and  $\eta\in C_{0}^2(G)$.
Then there exists $c_5>0$ such that
\begin{equation}\label{II-15}
 |(-\Delta)^s\eta(x)|\leq   \frac {c_5\|\eta\|_{C^2}}{1+|x|^{N+2s}} \qquad \forall\, x\in \R^N.
\end{equation}
Moreover, assume that $\eta\geq 0$   in $G$,
then
$(-\Delta)^s\eta\leq  0$ in $G^c$ and for any $\gd>0$
there exists $c_\gd>1$ independent of $\eta$ such that
\begin{equation}\label{II-16'}
\frac {\norm\eta_{L^1}}{c_\gd(1+|x|^{N+2s})}\leq   -(-\Delta)^s\eta(x)\leq  \frac {c_\gd\norm\eta_{L^1}}{1+|x|^{N+2s}},
\end{equation}
for $x\in\{z\in\R^N: \dist(z,G)\geq\gd\}.$
\es
\Proof Let $x\in G^c$ and $y\in \R^N$, then $\eta(x)-\eta(y)\leq  0$ and
hence $(-\Delta)^s\eta \leq  0$ in $  G^c$. For $y\in G$ and $x\in G^c$ satisfying $\dist(x,G)>\gd$, there exists $c_6>1$ such that
$$c_6^{-1}(1+|x|^{N+2s})\leq  |x-y|^{N+2s}\leq  c_6( 1+|x|^{N+2s}).$$
Together with
$$
 (-\Delta)^s\eta(x) = -a_{N,s}\int_{G}\frac{\eta(y)}{|y-x|^{N+2s}}dy \qquad\forall\, x\in G^c,
$$
one obtains the claim.\qeda\medskip

\nind\Remark Estimate (\ref{II-15})  has essentially been already obtained in \cite[Lemma 2.1]{BV} but we kept it for the sake of completeness.  (\ref{II-15}). Estimate (\ref{II-16}) is new and will be useful in the sequel. \medskip

\nind{\bf Proof of Theorem B}. Let $\gr>\gr'>0$ and $\gz\in C^2_0(B_{\gr}(z))$ such that $0\leq \gz\leq 1$ and $\gz =1$ on $B_{\gr'}(z))$. Then there holds
$$\BA {lll}
\myint{B_\gr(z)}{}u(t,x)\gz(x)dx=\myint{B_\gr(z)}{}u(T,x)\gz(x)dx+\myint{t}{T}\myint{B_\gr(z)}{}f(s,x,u)\gz(x)dx ds\\[4mm]
\phantom{-------------------}
+\myint{t}{T}\myint{\BBR^N}{}u(s,x)(-\Gd)^{s}\gz(x)dx ds.
\EA
$$
The function $\gz$ satisfies
$$\abs{(-\Gd)^{s}\gz(x)}\leq \myfrac{c_5\norm\gz_{C^2}}{1+\abs x^{N+2s}},\quad\forall\, x\in\BBR^N.
$$
Since $(t,x)\mapsto(1+\abs{x}^{N+2s})^{-1}u((t,x)\in L^1(Q_T)$, we infer that
\bel{II-15'}\BA {lll}\displaystyle
\lim_{t\to 0}\myint{B_\gr(z)}{}u(t,x)\gz(x)dx=+\infty,
\EA\ee
which implies the claim.
\qeda
\subsection{Pointwise estimates}

\nind{\bf Proof of Theorem C}. In what follows we characterize the singular set of the initial trace when the absorption reaction is subcritical, that is it satisfies (\ref{In-10}),  (\ref{In-11}) and (\ref{In-12}) hold. Under these two last assumptions for any bounded Radon measure in $\BBR^N$, it is proved in  \cite[Th 1.1]{CVW} that there exists a unique weak solution $u:=u_\gm$ to
\bel{II-16}\BA {lll}\displaystyle
\prt_tu+(-\Gd)^{s}u+t^\gb g(u)=0\quad&{\rm in }\;\ Q_\infty\\[1mm]
\phantom{--\prt_\Gd^{s}u+t^\gb g(u)}
\!u(0,\cdot)=\gm\quad&{\rm in }\;\ \R^N.
\EA\ee
We recall by a weak solution, we mean a function  $u\in L^1(Q_T)$  such that $t^\gb g(u)\in L^1(Q_T)$ for all $T>0$ satisfying
\bel{II-16+}\BA {lll}\displaystyle
\myint{0}{T}\myint{\BBR^N}{}\left[\left(-\prt_t\xi+(-\Gd)^s\xi\right)u+t^\gb g(u)\xi\right]dxdt=\myint{\BBR^N}{}\xi(0,x)d\gm(x)\qquad \forall\,\xi\in \BBY_{s,T },
\EA\ee
where $\BBY_{s,T }$ is the space of functions $\xi$ defined in $Q_\infty$ satisfying
$$
\BA{lll}
(i)\quad
\norm{\xi}_{L^1(Q_T)}+\norm{\xi}_{L^\infty(Q_T)}+\norm{\partial_t\xi}_{L^\infty(Q_T)}+\norm{(-\Delta)^s \xi}_{L^\infty(Q_T)}<+\infty,\\[2mm]

(ii)\quad \xi(T)=0\;{\rm and\;  for\; }
\;0<t<T,{ \rm there\; exist\; } M>0\;{ \rm and }\;\epsilon_0>0\; { \rm such\; that\; for\; } 0<\ge\leq \epsilon_0,\quad\\[1mm]
 \qquad\norm{(-\Delta)_\epsilon^s \xi(t,\cdot)}_{L^\infty(\R^N)}\leq M.
\EA
$$
Furthermore, if $\gm_j$ converges to $\gm$ weakly in the sense of measures, then $u_{\gm_j}$ converges to $u_\gm$ locally uniformly in $Q_\infty$. Up to translation we can assume that $z=0$. Since (\ref{In-13sup}) holds, for any $k>0$ there exist two sequences $\{t_n\}$ and $\{\gr_n\}$ converging to $0$ such that
 \bel{II-17}\BA {lll}\displaystyle
\myint{B_{\gr_n}}{}u(t_n,x)dx=k.
\EA\ee
\nind{\it Case 1: $\gb\leq 0$}. For $R>0$, let $v_n^R$ be the solution of
\bel{II-18}\BA {lll}\displaystyle
\prt_tv+(-\Gd)^{s}v+t^\gb g(v)=0\quad&{\rm in }\;\;(t_n,\infty)\ti B_R\\[1mm]
\phantom{--\prt_\Gd^{s}+t^\gb g(u)}
v(t,x)=0\quad&{\rm in }\;\;(t_n,\infty)\ti B^c_R\\[1mm]
\phantom{--\prt_\Gd^{s}+t^\gb g(u)}
\!v(t_n,.)=u(t_n,x)\chi_{_{B_{\gr_n}}}\quad&{\rm in }\;\;B_R,
\EA\ee
where $B_R$ denote the ball in $\R^N$ centered at origin with the radius $R$. 
By the comparison principle,  $u\geq v_n^R$ in $[t_n,\infty)\ti B_R$. We set $v_n^R(t,x)= v_n^R(\gt+t_n,x)=\tilde v_n^R(\gt,x)$. Since $\gb\leq 0$,  there holds
$$\prt_t\tilde v_n^R+(-\Gd)^{s}\tilde v_n^R+\gt^\gb g(\tilde v_n^R)\geq 0\quad{\rm in }\;\;(0,\infty)\ti B_R.
$$
Hence $\tilde v_n^R\geq u_n^R$ where $u_n^R$ is the solution of 
\bel{II-18*}\BA {lll}\displaystyle
\prt_tv+(-\Gd)^{s}v+t^\gb g(v)=0\quad&{\rm in }\;\;(0,\infty)\ti B_R\\[1mm]
\phantom{--\prt_\Gd^{s}+t^\gb g(v)}
v(t,x)=0\quad&{\rm in }\;\;(0,\infty)\ti B^c_R\\[1mm]
\phantom{--\!\prt_\Gd^{s}+t_n^\gb g(u_n)}
\!v(0,.)=u(t_n,x)\chi_{_{B_{\gr_n}}}\quad&{\rm in }\;\;B_R.
\EA\ee
Letting $R\to\infty$ we infer that $u_n^R$ increases and converges to the solution
$u_n^\infty$ of
\bel{II-19}\BA {lll}\displaystyle
\prt_tv+(-\Gd)^{s}v+t^\gb g(v)=0\qquad&{\rm in }\;\;(0,\infty)\ti \BBR^N\\[1mm]
\phantom{--\prt_\Gd^{s}+t^\gb g(v'')}
\!v(0,.)=u(t_n,x)\chi_{_{B_{\gr_n}}}\quad&{\rm in }\;\;\BBR^N,
\EA\ee
and there holds $u(t_n+\gt,x)\geq u_n^\infty(\gt,x)$ in $(0,\infty)\ti \BBR^N$. Letting $n\to\infty$ and using the above mentioned stability result, we obtain that $u_n^\infty$ converges to $u_{k\gd_0}$ and $u\geq u_{k\gd_0}$. Since it holds true for any $k$, the claim follows.\smallskip

\nind{\it Case 2: $\gb> 0$}. Clearly $u\geq v_n^\infty$ where $v_n^\infty$ satisfies 
\bel{II-18-2}\BA {lll}\displaystyle
\prt_tv+(-\Gd)^{s}v+t^\gb g(v)=0\qquad&{\rm in }\;\;(t_n,\infty)\ti \BBR^N\\[1mm]
\phantom{--\prt_\Gd^{s}+t^\gb g(v')}
\!v(t_n,.)=u(t_n,x)\chi_{_{B_{\gr_n}}}\quad&{\rm in }\;\;\BBR^N,
\EA\ee
Then $v_n^\infty(t,x)\leq \BBH_s[u(t_n,x)\chi_{_{B_{\gr_n}}}](t-t_n,x)$. Since $g$ satisfies $(\ref{In-11^*})$ it follows from \cite[Proof of Th 1.1]{CVW} that the set of function $\left\{g\left(\BBH_s[u(t_n,.)\chi_{_{B_{\gr_n}}}](.-t_n,.)\right)\right\}$ is uniformly integrable in $(t_n,\infty)\ti \BBR^N$ and it is the same with $\left\{g\left(v_n^\infty\right)\right\}$. Therefore, for any $T>0$, $\left\{t^\gb g\left(v_n^\infty\right)\right\}$
is uniformly integrable in $(t_n,T)\ti \BBR^N$. Hence $\{v_n^\infty\}$ converges locally in $(0,\infty)\ti \BBR^N$ to $u_{k\gd_0}$ and 
$u\geq u_{k\gd_0}$ as above.

$\phantom{--}$  \hfill$\Box$\medskip

\nind\Remark We will see in Theorem F that if $g(u)=u^p$ the concentration result holds under the mere condition  (\ref{In-11}) whatever is the sign of $\gb$. The difficulty in the case $\gb>0$ comes from the fact that the ball $B_{\gr_n}$ may shrink very quickly with $t_n$ and that a pointwise isolated singularity at $(\gt,z)$  with $\gt>0$ can be removable for equation (\ref{II-16}). In the power case we can control the rate of shrinking thanks to a Harnack-type inequality.\medskip

\nind{\bf Proof of Theorem D}. (i) {\it Proof of (\ref{In-13sup})}. Let $\gg\in C^2(\BBR)$ be a convex nondecreasing function vanishing on $(-\infty,0]$ such that $\gg(r)\leq r_+$. For $\ge>0$, let $U_\ge$ be the solution of
\bel{II-19-1}\BA {lll}\displaystyle
\prt_tU+t^\gb g(U)=0\quad&{\rm in }\;\;(\ge,+\infty)\\[1mm]
\phantom{-\prt\Gd\, g(u)}
\!\!U(\ge)=+\infty.
\EA\ee
Indeed, $$U_\ge(t)=\mathcal{G}^{-1}(\myfrac{t^{\gb+1}-\ge^{\gb+1}}{\gb+1}),$$
where $\mathcal{G}^{-1}$ is the inverse function of $\mathcal{G}$, see (\ref{In-14-2}).  
Then there holds
$$\BA {lll}\displaystyle
(-\Gd)^s\gg(u(t,.)-U_\ge(t))(x)=\gg'((u(t,x)-U_\ge(t))(-\Gd)^s u(x)\\[4mm]\phantom{-(-\Gd)^s\gg(u(t,.)-U_\ge(t))(x)}
-a_{N,s}\myfrac{\gg''(u(t,z_x)-U_\ge(t))}{2}\myint{\BBR^N}{}\myfrac{(u(t,y)-u(t,x))^2}{\abs{x-y}^{N+2s}}dy.
\EA$$
Notice that the integral is convergent if $t>\ge$,  since $\gg(u(t,\cdot)-U(t-\ge))=\gg(u(t,.)-U_\ge(t))$, where $0\leq u(t,\cdot)\leq U_\ge(t)$ and $u$ satisfies
$$\myint{\BBR^N}{}\myfrac{u(\cdot,x) dx}{1+|x|^{N+2s}}<+\infty\qquad{\rm a.e. \; in} \;\ (0,T).
$$
Then
$$\BA {lll}\displaystyle
\prt_t\gg(u(t,x)-U_\ge(t))+(-\Gd)^s\gg(u(t,.)-U_\ge(t))(x)\\[4mm]\phantom{------------}
\leq \gg'(u(t,x)-U_\ge(t))\cdot\left(\prt_tu(t,x)-\prt_tU_\ge(t)+(-\Gd)^s u(x)\right)\\[4mm]\phantom{------------}
\leq \gg'(u(t,x)-U_\ge(t))\cdot \left(t^\gb g(U_\ge(t))-f(t,x,u(t,x))\right)
\\[4mm]\phantom{------------}
\leq 0.
\EA$$
Therefore, $\gg(u(\cdot,\cdot)-U_\ge(\cdot))$ is a subsolution. Let $\eta\in C^\infty_0(\BBR^N)$ be a nonnegative function. Using \rlemma{test} we have that
$$
\left|\myint{\BBR^N}{}\gg(u(t,x)-U_\ge(t))(-\Gd)^s\eta dx\right|\leq c_5\norm\eta_{C^2}\myint{\BBR^N}{}\myfrac{u(t,x) dx}{1+\abs x^{N+2s}}.
$$
Since $u\in L_{loc}^1(0,T;\CL^s(\BBR^N))$, for almost all $s,t$ such that $\ge<s<t$, there holds
$$\BA {lll}\myint{\BBR^N}{}\gg(u(t,x)-U_\ge(t))\eta(x)dx+\myint{s}{t}\myint{\BBR^N}{}\gg(u(\gt,x)-U_\ge(t))(-\Gd)^s\eta(x) dxd\gt \\[4mm]
\phantom{\myint{\BBR^N}{}\gg(u(t,x)-U_\ge(t))\eta(x)dx-------}\leq \myint{\BBR^N}{}\gg(u(s,x)-U_\ge(s))\eta (x)dx.
\EA$$
Since $\gg(u(s,x)-U_\ge(s))\eta (x)\leq u(s,x)\eta(x)$ and $u(s,.)\eta\in L^1(\BBR^N)$, we get from the dominated convergence theorem that
 $$\displaystyle\lim_{s\downarrow\ge}\myint{\BBR^N}{}\gg(u(s,x)-U_\ge(s))\eta (x)dx=0.$$
Hence, letting $s\to\ge$ and  $\gg(r)\uparrow r_+$, we get
 \bel{II-20}\BA {lll}\displaystyle
\myint{\BBR^N}{}(u(t,x)-U_\ge(t))_+\eta(x)dx\leq \left|\myint{\ge}{t}\myint{\BBR^N}{}(u(\gt,x)-U_\ge(t))_+(-\Gd)^s\eta(x) dxd\gt \right|\\[4mm]
\phantom{\myint{\BBR^N}{}(u(t,x)-U_\ge(t))_+\eta(x)dx}
\leq c_5\norm{\eta}_{C^2}\myint{\ge}{t}\myint{\BBR^N}{}\myfrac{(u(\gt,x)-U_\ge(t))_+}{1+\abs x^{N+2s}} dxd\gt.
\EA\ee
 Next, for $n\geq 1$, we replace $\eta$ by $\eta_n(x)=\eta(n^{-1}x)$, where $0\leq\eta\leq 1$, $\eta(x)=1$ on $B_1$ and supp$(\eta)\subset B_2$. We can also assume that $\eta$ is radially decreasing and $\eta(0)=1$. Since   $\norm{\eta_n}_{C^2}\leq \norm{\eta}_{C^2}$, we obtain from
 (\ref{II-20}) and the monotone convergence theorem that the following holds for almost all $t\in (\ge,T)$
  \bel{II-21}\BA {lll}\displaystyle
\myint{\BBR^N}{}(u(t,x)-U_\ge(t))_+dx\leq c_5\norm{\eta}_{C^2}\myint{\ge}{t}\myint{\BBR^N}{}\myfrac{(u(\gt,x)-U_\ge(\gt))_+}{1+\abs x^{N+2s}} dxd\gt.
\EA\ee
 This inequality implies that $(u(t,.)-U_\ge(t))_+\in L^1(\BBR^N)$ for almost all $t\in (\ge,T)$. We set
 $$\Psi_\ge(t)=\myint{\ge}{t}\myint{\BBR^N}{}\myfrac{(u(\gt,x)-U_\ge(\gt))_+}{1+\abs x^{N+2s}} dxd\gt.
 $$
 Then
  $$\Psi'_\ge(t)=\myint{\BBR^N}{}\myfrac{(u(t,x)-U_\ge(t))_+}{1+\abs x^{N+2s}} dx\leq \myint{\BBR^N}{}(u(t,x)-U_\ge(t))_+dx\leq
  c_5\norm{\eta}_{C^2}\Psi_\ge(t).
 $$
 Since $\Psi_\ge(\ge)=0$ we obtain $\Psi_\ge(t)=0$ on $(0,T)$, hence $u(t,x)\leq U_\ge(t)$ a.e. on $(\ge,T)\ti\BBR^N$. Letting $\ge\to 0$, we get the claim.\smallskip

 \nind (ii) {\it End of the proof}. Because of Theorem B it is sufficient to prove that if (\ref{In-14-4}) holds true, then $U\in L^1(0,1)$. Indeed, we recall that
 $$ \mathcal{G}(s)=\myint{s}{\infty}\myfrac{d\tau}{g(\tau)}.
 $$
Clearly, $\mathcal{G}$ is an decreasing diffeomorphism from $\BBR^*_+$ onto $(0,\Gf(0))$ and $U(t)=\mathcal{G}^{-1}\left(\frac{t^{\gb+1}}{\gb+1}\right)$.   Set $U(t)=s$, then $t=\left((\gb+1)\mathcal{G}(s)\right)^{\frac{1}{\gb+1}}$ and we get
 $$\BA {lll}
 \myint{0}{1}U(t)dt=\myint{\infty}{U(1)}s\mathcal{G}'(s)\left((\gb+1)\mathcal{G}(s)\right)^{-\frac{\gb}{\gb+1}}ds\\[4mm]
 \phantom{\myint{0}{1}U(t)dt}
 =(\gb+1)^{-\frac{\gb}{\gb+1}}\myint{U(1)}{\infty}\myfrac{sds}{g(s)\left(\myint{s}{\infty}\myfrac{d\gt}{g(\gt)}\right)^{\frac{\gb}{\gb+1}}}<+\infty,
 \EA$$
 which completes the proof.\qeda\medskip

The following weight function plays an important role in the description of the initial trace problem for positive solutions of the fractional heat equation
  \bel{G-0}\BA {lll}\displaystyle
\Gf(x)=\myfrac{1}{\left(1+\left(\abs x^2-1\right)^4_+\right)^{\frac{N+2s}{8}}},\quad\forall\, x\in \BBR^N.
\EA\ee
It has the remarkable property that
  \bel{G-1}\BA {lll}\displaystyle
-c_6\Phi(x)\leq (-\Gd)^s\Phi(x)\leq c_6\Phi(x),\quad\forall\, x\in \BBR^N,
\EA\ee
for some constant $c_6>0$ (see \cite{BSV}, \cite{BV}). Furthermore, for some $c_7>1$,
  \bel{G-2}\BA {lll}\displaystyle
\myfrac{1}{c_7(1+\abs x^{N+2s})}\leq \Phi(x)\leq \myfrac{c_7}{1+\abs x^{N+2s}},\quad\forall\, x\in \BBR^N.
\EA\ee

\blemma{estim} Let  $f:\BBR^*_+\ti\BBR^N\ti\BBR_+\rightarrow\BBR_+$ be a Caratheodory function which satisfies (\ref{In-2}) and is nondecreasing with respect to  the variable $u$. For given $u_0\in L^1(\BBR^N)$ is nonnegative, problem
  \bel{G-3}\BA {lll}\displaystyle
\prt_tu+(-\Gd)^s u+f(t,x,u)=0&\quad{\rm in}\;\;Q_\infty,\\[1mm]
\phantom{--------\ \ \, }
u(0,\cdot)=u_0&\quad{\rm in}\;\;\BBR^N
\EA\ee
 has a unique weak solution $u\in C(\BBR^+;L^1(\BBR^N))$
satisfying that 
  \bel{G-4}\BA {lll}\displaystyle
\myint{\BBR^N}{}u(t,x)\Phi(x)dx+
\myint{0}{t}\myint{\BBR^N}{}\left(u(s,x)(-\Gd)^s\Phi(x)+f(s,x,u)\Phi(x)\right)dx ds=\myint{\BBR^N}{}u_0(x)\Phi(x)dx.
\EA\ee
\es
\Proof  Since $u$ is a weak solution of (\ref{G-3}) and the function $\Phi$ satisfies the assumptions (i)-(ii) in \cite[Def. 1.1]{CVW}, we get (\ref{G-4}).
\qeda

\bcor{estim-2} Assume that $f$ satisfies the assumptions of \rlemma{estim} and that inequalities (\ref{In-10})-(\ref{In-11})-(\ref{In-12}) hold. Then for any nonnegative measure $\gm$ in $\BBR^N$ verifying
  \bel{G-6}\BA {lll}\displaystyle
\myint{\BBR^N}{}\Phi(x)\,d\gm(x)<+\infty,
\EA\ee
there exists a weak solution $u\in C_b(\BBR^+;\CL^s(\BBR^N))\cap L^1(\BBR^+;\CL^s(\BBR^N))$  of (\ref{G-3}) in the sense that
  \bel{G-7}\BA {lll}\displaystyle
\myint{0}{t}\myint{\BBR^N}{}\left[-(\prt_t\xi+(-\Gd)^s\xi)u+\xi f(s,x,u)\right]dxds+\myint{\BBR^N}{}u(t,x)\xi (t,x)dx=\myint{\BBR^N}{}\xi (0,x)d\gm(x),       \EA\ee
for any $\xi\in C^2_{0}(\overline Q_T)$ satisfying the assumptions (i)-(ii) in \cite[Def. 1.1]{CVW}. Furthermore,
  \bel{G-8}\BA {lll}\displaystyle
\myint{\BBR^N}{}u(t,x)\Phi(x)dx+
\myint{0}{t}\myint{\BBR^N}{}\left[u(s,x)(-\Gd)^s\Phi(x)+f(s,x,u)\Phi(x)\right]dx ds
=\myint{\BBR^N}{}\Phi(x)d\gm(x).
\EA\ee
\es
\Proof By the assumptions on $f$ and for any $n>0$, it follows from \cite[Th. 1.1]{CVW} that
  \bel{G-8'}\BA {lll}\displaystyle
\prt_tu+(-\Gd)^s u+f(t,x,u)=0&\quad{\rm in}\;\;Q_\infty\\[1mm]
\phantom{---------}
u(0,\cdot)=\gm_n:=\chi_{B_n}u_0&\quad{\rm in}\;\;\BBR^N,
\EA\ee
has a unique solution $u_n\in L^1(Q_T)$ verifying $f(\cdot,\cdot,u_n)\in L^1(Q_T)$.
If $\gr_{k}$ is a sequence of mollifiers with compact supports and $\gm_{n,k}=(\chi_{B_n}u_0)\ast\gr_k$, the sequence $\{u_{n,k}\}$ of weak solutions of
  \bel{G-9}\BA {lll}\displaystyle
\prt_tu+(-\Gd)^s u+f(t,x,u)=0&\quad{\rm in}\;\;Q_\infty\\[1mm]
\phantom{--------\ \ \ }
u(0,\cdot)=\gm_{n,k}&\quad{\rm in}\;\;\BBR^N,
\EA\ee
then $u_{n,k}$ satisfies that
  \bel{G-10}\BA {lll}\displaystyle
  \myint{\BBR^N}{}u_{n,k}(t,x)\Phi(x)dx+
\myint{0}{t}\myint{\BBR^N}{}\left[u_{n,k}(s,x)(-\Gd)^s\Phi(x)+f(s,x,u_{n,k})\Phi(x)\right]dx ds\\[3mm]
\phantom{---\myint{\BBR^N}{}\left[\myfrac{c_{N,s}u_{n,k}(s,x)}{(1+\abs x^2)^{N+2s}}+\myfrac{f(s,x,u_{n,k})}{(1+\abs x^2)^{N-2s}}\right]dx ds}
=\myint{\BBR^N}{}\gm_{n,k}(x)\Phi(x)dx.
\EA\ee
When $k\to\infty$, we know from the proof of \cite[Th. 1.1]{CVW} that, up to a subsequence, $\{u_{n,k}\}_k$ converges a.e. in $Q_T$ to some function $u_n$, $\{f(\cdot,\cdot,u_{n,k})\}_k$ converges a.e. to $\{f(\cdot,\cdot,u_{n})\}$ and that $\{u_{n,k}\}_k$ and $\{f(\cdot,\cdot,u_{n,k})\}_k$ are uniformly integrable in $L^1(Q_T)$. Furthermore
$u_n\in C([0,T];L^1(\BBR^N))$ and for any $t\in (0,T]$, $\{u_{n,k}(t,\cdot)\}_k$ converges to $u_{n}(t,\cdot)$ in $L^1(\BBR^N)$. This implies that
  \bel{G-10'}\BA {lll}\displaystyle
 \myint{\BBR^N}{}u_{n}(t,x)\Phi(x)dx+
\myint{0}{t}\myint{\BBR^N}{}\left[u_{n}(s,x)(-\Gd)^s\Phi(x)+f(s,x,u_{n})\Phi(x)\right]dx ds=\myint{\BBR^N}{}\Phi(x)d\gm_n(x).
\EA\ee
Furthermore,
\bel{G-11}\BA {lll}\displaystyle
\myint{0}{t}\myint{\BBR^N}{}\left[-(\prt_t\xi+(-\Gd)^s\xi)u_n+\xi f(s,x,u_n)\right]dxds+\myint{\BBR^N}{}u_n(t,x)\xi (t,x)dx=\myint{\BBR^N}{}\xi (0,x)d\gm_n(x),      
 \EA\ee
for all $\xi\in C^2_{0}(\overline Q_T)$ satisfying the assumptions (i)-(ii) in \cite[Def. 1.1]{CVW}.
When $n\to\infty$, $u_n\uparrow u$ and  $f(s,x,u_{n})\uparrow f(s,x,u)$. Using the monotone convergent theorem we see that
$u$ satisfies (\ref{G-8}),and that the sequences $\{u_n\}_n$ and $\{f(\cdot,\cdot,u_{n})\}_n$ converges to $u$ and $f(\cdot,\cdot,u)$ in $L^1(0,T;\CL^s(\BBR^N))$  respectively. Using estimate (\ref{II-15}) we can let $n$ to infinity in  (\ref{G-11}) and obtain (\ref{G-7}).\qeda\medskip

As it is pointed out in \cite{BSV}, the weight function $\Gf$ plays a role similar to an eigenfunction of $(-\Gd)^s$. We prove a backward-forward uniqueness result for solutions of (\ref{In-1}) inspired from \cite[Lemma 4.2]{BSV}.

\bth{unic} Assume that $u\mapsto f(t,x,u)$ is locally Lipschitz continuous on $\BBR$, uniformly with respect to $x\in\BBR^N$ and locally uniformly with respect to $t\in \BBR^*_+$. If $u_1$ and $u_2$ belong to $L^1_{loc}(\BBR^*_+;\CL^s(\BBR^N))\cap L^\infty_{loc}(\BBR^*_+;L^\infty(\BBR^N))$ and are weak solutions of (\ref{In-1}) in $Q_T$
which coincide for $t=t_0>0$, then $u_1=u_2$ in $Q_T$.
\es
\Proof For any $0<\ge<t_0<T<\infty$, $u_1$ and $u_2$ are uniformly bounded in $[\ge,T]\ti\BBR^N$. Hence the function $D$ defined by
$$D(t,x)=\left\{\BA{lll}\myfrac{f(t,x,u_1(t,x))-f(t,x,u_2(t,x))}{u_1(t,x)-u_2(t,x)}\quad&{\rm if }\;\;u_1(t,x)\neq u_2(t,x)\\[2mm]
0&{\rm if }\;\;u_1(t,x)= u_2(t,x)
\EA\right.
$$
is bounded in $[\ge,T]\ti\BBR^N$ by some constant $M=M(\ge,T)>0$. Set $w=u_1-u_2$, it satisfies
$$\prt_tw+(-\Gd)^s w+D w=0 \quad{\rm in }\;\;Q_T,
$$
and is uniformly bounded in $[\ge,T]\ti\BBR^N$. Hence
$$\myfrac{d}{dt}\myint{\BBR^N}{}w(t,x)\Gf(x)dx+\myint{\BBR^N}{}w(t,x)(-\Gd)^s\Gf(x)dx+\myint{\BBR^N}{}D(t,x)w(t,x)\Gf(x)dx=0.
$$
Using (\ref{G-1}) we get
\bel{G-12}\BA {lll}\displaystyle
-(c_5+M)\myint{\BBR^N}{}w(t,x)\Gf(x)dx  \leq \myfrac{d}{dt}\myint{\BBR^N}{}w(t,x)\Gf(x)dx\leq (c_5+M)\myint{\BBR^N}{}w(t,x)\Gf(x)dx.
\EA\ee
This implies
\bel{G-13}\BA {lll}\displaystyle
-(c_5+M)\myint{\BBR^N}{}w(t,x)\Gf(x)dx  \leq \myfrac{d}{dt}\myint{\BBR^N}{}w(t,x)\Gf(x)dx\leq (c_5+M)\myint{\BBR^N}{}w(t,x)\Gf(x)dx,
\EA\ee
and
\bel{G-14}\BA {lll}\displaystyle
(i)\qquad\myint{\BBR^N}{}w(t,x)\Gf(x)dx\leq e^{(c_5+M)(t-s)}\myint{\BBR^N}{}w(s,x)\Gf(x)dx,\\[3mm]
(ii)\qquad e^{(c_5+M)(s-t)}\myint{\BBR^N}{}w(s,x)\Gf(x)dx\leq \myint{\BBR^N}{}w(t,x)\Gf(x)dx,
\EA\ee
for all $\ge\leq s\leq t\leq T$. Taking $s=t_0$ in (i) and $t=t_0$ in (ii) yields $w\equiv 0$ in $[\ge,T]\ti\BBR^N$.\qeda\medskip

\nind{\bf Proof of Theorem E}. By Theorem D we know that $u\leq U$. If there exists some $(t_0,x_0)\in Q_T$ such that $u((t_0,x_0))= U(t_0)$, then
either $u((t_0,x))= U(t_0)$ for all $x\in\BBR^N$, or
$$(-\Gd)^s(u-U)(t_0,x_0)<0\qquad\forall\, x\in\BBR^N.
$$
Since $f(t,x,u)-t^\gb g(U)\geq 0$ and $\prt_t(u-U)(t_0,x_0)=0$ we infer that $u((t_0,.))\equiv U(t_0)$. Since $g$ is nondecreasing this situation is impossible, hence $u((t_0,.))= U(t_0)$. Since $g$ is locally Lipschitz continuous, this implies $u=U$ in $Q_T$ by \rth{unic}.
\qeda\medskip

A straightforward consequence of Theorems B, C and D is the next statement.

\bcor{B3-4} Let $f(t,x,r)=t^\gb g(r)$, where $\gb>-1$ and $g:\BBR_+\rightarrow\BBR_+$ is continuous and nondecreasing and satisfies
(\ref{In-11}), (\ref{In-14-1}) and (\ref{In-14-4}). If $u$ is a nonnegative of (\ref{In-1}) in $Q_T$ belonging to $L^1_{loc}(0,T;\CL^s(\BBR^N))$ such that $\CS_u\neq\emptyset$, there holds
\bel{G-15}\BA {lll}u(t,x)\geq u_{\infty,z}(t,x)=u_{\infty,0}(x-z,t)\qquad\forall\, (t,x)\in Q_T.
\EA\ee
\es

\section{The case $f(t,x,u)=t^\gb u^p$}

We denote by $(-\Gd)_{_{\BBR^\gk}}^s$ the fractional Laplacian in $\BBR^\gk$ and $(-\Gd)_{_{\BBR^N}}^s=\left(-\Gd\right)^{s}$. The following standard lemma will be useful in the sequel.
\blemma{sum} Let $1\leq \gk\leq N-1$ be an integer. If $u\in C^2(\BBR^\gk)\cap\CL^{s}(\BBR^\gk)$ and $\tilde u(x_1,x')= u(x_1)$ for $(x_1,x')\in\BBR^\gk\ti\BBR^{N-\gk}$, then
\bel{C-1+}\BA {lll}\displaystyle
(-\Gd)^s \tilde u(x_1,x')=(-\Gd)_{_{\BBR^\gk}}^su(x_1).
\EA\ee
\es
\Proof This more or less well known lemma is based upon the explicit value of the constant $a_{N,s}$ in the definition of $(-\Gd)^s$. For the sake of completeness we give here the proof.
$$\BA {lll}(-\Gd)^s \tilde u(x_1,x')=a_{N,s}\myint{\BBR^\gk}{}\myint{\BBR^{N-\gk}}{}\myfrac{u(x_1)-u(y_1)}{\left((x_1-y_1)^2+\abs{x'-y'}^2\right)^{\frac N2+s}}dy'dy_1\\[4mm]\phantom{(-\Gd)^s \tilde u(x_1,x')}
=a_{N,s}\myint{\BBR^\gk}{}\left(\myint{\BBR^{N-\gk}}{}\myfrac{dy'}{\left((x_1-y_1)^2+\abs{y'}^2\right)^{\frac N2+s}}\right)\left(u(x_1)-u(y_1)\right)dy_1
\\[4mm]\phantom{(-\Gd)^s \tilde u(x_N,x')}
=a_{N,s}\left(\myint{\BBR^{N-\gk}}{}\myfrac{dz'}{\left(1+\abs{z'}^2\right)^{\frac N2+s}}\right)\myint{\BBR^\gk}{}\myfrac{u(x_1)-u(y_1)}{\abs{x_1-y_1}^{\gk+2s}}dy_1
\\[4mm]\phantom{(-\Gd)^s \tilde u(x_N,x')}
=\myfrac{a_{N,s}}{a_{\gk,s}}\left(\left|S^{N-1-\gk}\right|\myint{0}{\infty}\myfrac{r^{N-\gk-1}dr}{\left(1+r^2\right)^{\frac N2+s}}\right)
(-\Gd)_{_{\BBR^\gk}}^su(x_1).
\EA$$
Since
$$\left|S^{N-1-\gk}\right|=\myfrac{2\gp^{\frac{N-1-\gk}{2}}}{\Gg(\frac{N-1-\gk}{2})},
$$
and (see e.g. \cite[p. 103]{Spie})
$$\myint{0}{\infty}\myfrac{r^{N-\gk-1}dr}{\left(1+r^2\right)^{\frac N2+s}}=\myfrac{1}{2}B\left(\frac{N-\gk}{2},\frac{\gk}{2}+s\right)
=\myfrac{1}{2}\myfrac{\Gg(\frac \gk2+s)\Gg(\frac {N-\gk}2)}{\Gg(\frac N2+s)},
$$
by Euler's formula, where $B$ denotes beta function, we deduce that
$$\myfrac{a_{\gk,s}}{a_{N,s}}=\left|S^{N-1-\gk}\right|\myint{0}{\infty}\myfrac{r^{N-\gk-1}dr}{\left(1+r^2\right)^{\frac N2+s}},
$$
which yields (\ref{C-1+}).
\qeda\medskip

The next statement is a straightforward consequence.

\bcor{sum+} Assume that $\displaystyle u(x)=u(x_1,...x_N)=\sum_{j=1}^{N}u_j(x_j)$, then
\bel{C-3}\displaystyle
(-\Gd)^s \tilde u(x)=\sum_{j=1}^{N}\left((-\Gd)_{_{\BBR}}^su_j\right)(x_j).
\ee
\es

\subsection{Proof of Theorem F}

By Theorem E there holds $u(t,x)\leq ct^{-\frac{1+\gb}{p-1}}$ for some $c_*=c_*(\gb,p)>0$, hence $u$ satisfies 
\bel{F-1} \prt_tu+(-\Gd)^s u +\myfrac{c(t,x)}{t}u=0,
\ee
where $0\leq c(t,x)=t^\gb u^{p-1}(t,x)\leq c_*^{p-1}:=C_*$. Let $d_s$ be the fractional parabolic pseudo-distance (i.e. the triangle inequality holds up to a multiplicative constant if $s<\frac 12$) in $\BBR^N\ti\BBR$,  
$$d_s((t,x),(s,y))=\sqrt{|x-y|^2+\abs{t-s}^{\frac{1}{s}}}.
$$

\blemma {LF-1}If $z\in\CS_u$, there holds 
\bel{F-2} \limsup_{d_s((t,x),(0,z))\to 0} t^{\frac{N}{2s}}u(t,x)=\infty.
\ee
\es
\Proof If (\ref{F-2}) does not hold there exists $m,\ge_0>0$ such that 
$$u(t,x)\leq mt^{-\frac{N}{2s}}\qquad\forall (t,x)\;{\rm s.t.}\;\,|x-z|^2+t^{\frac{1}{s}}\leq \ge^2_0.
$$
Hence 
$$\gg(t,x):=t^\gb u^{p-1}(t,x)\leq m^{p}t^{\gb-\frac{N (p-1)}{2s}}\qquad\forall (t,x)\in (0,t_1)\ti B_{\ge_1}(z),
$$
where $\ge_1=\frac{\ge_0}{2}$ and $t_1=\left(\frac 34\right)^{2s}\ge_0^s$. By assumption $p<p^*_\gb$, or equivalently 
$\gb-\frac{N (p-1)}{2s}>-1$. Hence $\tilde\gg(t):=\norm{\gg(.,t)}_{L^\infty(B_{\ge_1}(z))} \in L^1(0,t_1)$.  We write the equation satisfied by $u$ 
in $B_{\ge_1}(z)\ti (0,t_1)$ in the form
\bel{F-3} \prt_tu+(-\Gd)^s u +\gg(t,x) u=0,
\ee
and, as in the proof of Theorem A, we take for test function $\phi^q$, where $q\geq 2$ and $\phi_\ge=\phi_{B_{\ge}(z)}$ is the first normalized 
eigenfunction of $(-\Gd)^s$ in $H^s_0(\phi_{B_{\ge}(z)})$ for some $0<\ge<\ge_1$. If $\gl_{\ge}$ is the corresponding eigenvalue, we obtain as in \rprop{loc1}, 
$$\myfrac{d}{dt}e^{q\gl_{\ge}t}\myint{B_{\ge}(z)}{}u\phi_\ge^q dx+e^{q\gl_{\ge}t}\tilde\gg(t)\myint{B_{\ge}(z)}{}u\phi_\ge^q dx\geq 0.
$$
If we put $X(t)=e^{q\gl_{\ge}t}\myint{B_{\ge}(z)}{}u\phi_\ge^q dx$, then $X'+\tilde\gg(t) X\geq 0$ on $(0,t_1)$, which implies that 
the function $t\mapsto e^{\int_0^{t}\tilde\gg(s)ds}X(t)$ is increasing on $(0,t_1)$. Hence 
$$\displaystyle\lim_{t\to 0}\myint{B_{\ge}(z)}{}u(.,t)\phi_\ge^q dx\leq e^{q\gl_{\ge}t_1+\int_0^{t_1}\tilde\gg(s)ds}\myint{B_{\ge}(z)}{}u(.,t_1)\phi_\ge^q dx,
$$
which implies that $z\in \CR_u$, contradiction.
\qeda\medskip

Notice that the above lemma contains a result which is interesting in itself.

\bcor{opt} If $\gg$ is a measurable function in $\BBR^N\ti (0,T)$ such that for any compact set $K\subset\BBR^N$ 
the function $\displaystyle\tilde\gg_K(t):={\rm ess}\sup_{\!\!\!\!\!x\in K}|\gg(t,x)|$ is integrable on $(0,T)$, then any nonnegative function $u\in L^{1}_{loc}(0,T;\CL^s(\BBR^N)$ satisfying (\ref{F-3}) admits an initial trace $\gn$ which is a nonegative Radon measure in $\BBR^N$.
\es

The next result is an Harnack-type inequality valid for positive solutions of (\ref{In-15}). For the mere fractional heat equation, two-sided Harnack inequalitis are proved in \cite{BV} and \cite{BSV}.

\blemma {LF-2} Let $\gth>0$ and $w$ be a nonnegative solution of  (\ref{In-15}) in $Q_\infty$. Then for any $t>\frac{t}{2}\geq s>\frac{t}{4}$ and
$x,y\in\BBR^N$ such that $\abs{x-y}\leq \gth t^{\frac{1}{2s}}$, there holds
\bel{F-4} \BA{lll}
w(t,y)\geq Mw(x,s),
\EA\ee
where $M>0$ depends on $N$, $s$, $\gb$, $p$ and $\gth$.
\es
\Proof Since $w$ satisfies 
\bel{F-5}
\prt_tw+(-\Gd)^s w+\myfrac{C_*}{t}w\geq 0,
\ee
$t^{C_*}w(t,.)$ is a supersolution of the fractional diffusion equation $\prt_tv+(-\Gd)^s v=0$, hence 
$$\BA{lll}t^{C_*}w(t,x)\geq \left(\myfrac{t}{4}\right)^{C_*}\myint{\BBR^N}{}H_s(\frac{3t}{4},x-z)w(\frac{t}{4},z)dz,
\EA$$
which implies, thanks to identity (\ref{cks}),
$$\BA{lll}w(t,x)\geq \myfrac{3t}{c_34^{C_*+1}}\myint{\BBR^N}{}\myfrac{w(\frac{t}{4},z)dz}{(\frac{3t}{4})^{1+\frac{N}{2s}}+\abs{x-z}^{N+2s}}.
\EA$$
Since $w$ is a subsolution of the fractional diffusion equation,
$$\BA{lll}w(s,y)\leq \myint{\BBR^N}{}H_s(s-\frac{t}{4},y-z)w(\frac{t}{4},z)dz\leq c_3(s-\frac{t}{4})
\myint{\BBR^N}{}\myfrac{w(\frac{t}{4},z)dz}{(s-\frac{t}{4})^{1+\frac{N}{2s}}+\abs{y-z}^{N+2s}}.
\EA$$
Hence
\bel{F-6}w(s,y)\leq \myfrac{c_3^24^{C_*+1}(s-\frac{t}{4})}{3t} \left(\sup_{z\in\BBR^N}\myfrac{(\frac{3t}{4})^{1+\frac{N}{2s}}+\abs{x-z}^{N+2s}}{(s-\frac{t}{4})^{1+\frac{N}{2s}}+\abs{y-z}^{N+2s}}\right)w(t,x)
\ee
If we assume that $|x-y|\leq \gth t^{\frac{1}{2s}}$ for some $\gth>0$, we obtain the claim.\qeda \medskip

\nind{\bf End of the proof of Theorem F}. By \rlemma{LF-1} there exists a sequence $\{t_n,x_n\}\subset Q_\infty$ converging to $(0,z)$ such that 
\bel{F-7*}
u(s_n,x_n)\geq ns_n^{-\frac{N}{2s}}.
\ee
By \rlemma {LF-2}, there holds with $t_n=2s_n$, 
\bel{F-7}
u(t_n,x)\geq Mnt_n^{-\frac{N}{2s}}\geq cn H_{s}(t_n,x-x_n)\qquad\forall x\in \BBR^N\;{\rm s.t.}\;\abs{x-x_n}\leq \gth t_n^{\frac{1}{2s}}.
\ee
for some $c>0$ depending on $M$ and $\gth$. This implies 
\bel{F-8}
\myint{\abs{x-x_n}\leq \gth t_n^{\frac{1}{2s}}}{}u(t_n,x)dx\geq cn\myint{\abs{x-x_n}\leq \gth t_n^{\frac{1}{2s}}}{}H_{s}(t_n,x-x_n)dx\geq c'n.
\ee
Then for any $k\in (0,n)$ there exists $k_n\in (0,n)$ such that 
\bel{F-9}
ck_n\myint{\abs{x-x_n}\leq \gth t_n^{\frac{1}{2s}}}{}H_{s}(t_n,x-x_n)dx=k.
\ee
Using $(\ref{In-24})$ 
$$\myint{\abs{x-x_n}\leq \gth t_n^{\frac{1}{2s}}}{}H_{s}(t_n,x-x_n)dx=\myint{B_\gth}{}\tilde H_s(y)dy,
$$
hence $k_n\leq ck$ for some $c>0$ independent of $n$.  
Set $\gu_n(x)=ck_nH_{s}(t_n,.-x_n)\chi_{B_{\gth t_n^{\frac{1}{2s}}}(x_n)}$. Then $u$ is bounded from below in $(t_n,\infty)\ti\BBR^N$ by the function $u_n$  which satisfies 
\begin{equation}\label{F-10}
 \arraycolsep=1pt
\begin{array}{lll}
\partial_t u_n + (-\Delta)^s  u_n+ t^\beta u_n^p=0\quad & {\rm in}\quad (t_n,\infty)\ti\BBR^N\\[1mm]
\phantom{------,,,\,\,}
u_n(t_n,.)=\gu_n,
\end{array}
\end{equation}
which in turn, satisfies
$$u_n(t,x)\leq \BBH_s[\gu_n](t-t_n,x)\leq ck_nH_{s}(t,x-x_n)\quad  {\rm in}\quad (t_n,\infty)\ti\BBR^N.$$
Since $p<p^*_\gb$, it is proved in \cite[Th 1.1]{CVW} that the set of functions $\{t^\gb \left(ck_nH_{s}(.,.-x_n)\right)^p\}$ is uniformly integrable in $Q_\infty$, and this property is shared by the set $\{t^\gb \left(u_n\right)^p\}$. Because $u_n(t_n,.)\to k\gd_z$ it follows that $u_n\to u_{k\gd_z}$ locally uniformly in $Q_\infty$, and $u\geq u_{k\gd_z}$. Since $k$ is arbitrary, the claim follows.\qeda

\subsection{Proof of Theorem G  (i)}

When $f(t,x,u)=t^\gb g(u):=t^\gb u^p$, conditions (\ref{In-14-1}) and (\ref{In-14-4}) are fulfilled when $p>1$ and $p>\gb+2$ respectively. Condition $1<p\leq p^{**}_\gb$ is not compatible with $p>\gb+2$, and condition $p^{**}_\gb<p<p^{*}_\gb$ necessitates $\gb+2<1+\frac{2s(1+\gb)}{N}$, equivalently $\gb+1<\frac{2s}{N}$.\medskip

 \nind{\it Step 1. The case $1<p<1+\frac{2s(1+\beta)}{N+2s})$}. Let $z\in \CS_u$. Since $r\mapsto r^p$ satisfies  (\ref{In-11}) there holds
 $u\geq u_{z,\infty}$ by Theorem C. Since $u_{z,\infty}=U_{p,\gb}$ by (\ref{In-17}), we obtaind that $u\geq U_{p,\gb}$. If we assume that
 $u\in L^1_{loc}(0,T;\CL^s(\BBR^N))$ then $u= U_{p,\gb}$ by Theorem D.
\medskip

 \nind{\it Step 2. The case $p=1+\frac{2s(1+\beta)}{N+2s}$}. We set $u_\infty=u_{0,\infty}$. From \cite[Lemma 5.3]{CVW} ,
\begin{equation}\label{3.5}
u_\infty(t,x)\geq  \frac{c_7 t^{-\frac{N+2s}{2s} }}{1+|t^{-\frac1{2s}}x|^{N+2s}} \qquad \forall\, (t,x)\in (0,1)\times\R^N,
\end{equation}
where $c_7>0$.
Since
$$\myint{B_1}{}\frac{t^{-\frac{N+2s}{2s} }dx}{1+|t^{-\frac1{2s}}x|^{N+2s}}=t^{-1}\myint{0}{1}\myfrac{r^{N-1}dr}{1+r^{N+2s}},
$$
it follows from (\ref{3.5}) that
\begin{equation}\label{3.6}
\int_0^1\int_{B_1(0)} u_\infty(t,x)dxdt=+\infty.
\end{equation}
By \rprop {sym} in Appendix,  $x\mapsto u_k(t,x)$ is  radially symmetric and decreasing, so is $u_\infty$. Therefore, if we are able to prove that there exists $x\in\R^N\setminus\{0\}$ such that $\displaystyle\lim_{t\to 0}u_\infty(t,x)=\infty$, it will imply
$$\displaystyle\lim_{t\to 0}u_\infty(t,z)=\infty\qquad{\rm uniformly\; with\;  respect\;  to\;  z\;  in }\;\;\overline B_{\abs x}(0).
$$
Hence $\overline B_{\abs x}(0)\subset \CS_{u_\infty}$ and by Theorem C,
\begin{equation}\label{C4}
u_\infty(t,x)\geq u_{z,\infty}(t,x)=u_{\infty}(t,x-z)\qquad\forall\, z\in \overline B_{\abs x}(0)
\end{equation}
Because $u_\infty$ is  radially symmetric and decreasing, it implies that
\begin{equation}\label{C5}
u_\infty(t,x)=u_{\infty}(t,x-z)\qquad\forall\, z\in \overline B_{\abs x}(0).
\end{equation}
By iterating this process we infer that $u_\infty(t,x)$ is indeed independent of $x$ and tends to $\infty$ when $t\to 0$. It coincides therefore to the maximal solution  $U_{p,\gb}$ of (\ref{In-14-3}) with $g(u)=u^p$.\smallskip

Henceforth we are lead to prove that $\CS_{u_\infty}\cap\BBR^N\setminus\{0\}\neq \emptyset$.  We proceed by contradiction in supposing that it does not hold, and let $x_0\in\CS^c_{u_\infty}\cap\BBR^N\setminus\{0\}$. Then $\displaystyle\limsup_{t\to 0}u_\infty(t,x_0)<\infty$ and
\begin{equation}\label{C6}\displaystyle \sup_{0<t\leq 1}u_\infty(t,x)\leq \displaystyle\sup_{0<t\leq 1}u_\infty(t,x_0):=M<\infty
\quad{\rm uniformly\; with\;  respect\;  to\;  {\it x}\;  in }\;\;\overline B^c_{\abs {x_0}}(0).
\end{equation}
By rescaling we can assume that $\abs {x_0}=1$. Let $\bar x\in \overline B^c_{3}$ and $\eta\in C^2_{0}(B_1(\bar x))$ such that $\eta\geq 0$ and $\eta=1$ on $\overline B_{\frac 12}(\bar x)$. We denote
$$X_1(t)=\myint{\BBR^N}{}\!\!\!u_\infty(t,x)\eta(x) dx,\;\; Y_1(t)=t^\gb\myint{\BBR^N}{}\!\!\!u^p_\infty(t,x)\eta(x) dx,\;\; Z_1(t)=\myint{\BBR^N}{}\!\!\!u_\infty(t,x)(-\Gd)^s\eta(x) dx,
$$
and there holds
\begin{equation}\label{C7}\displaystyle
X'_1(t)+Z_1(t)+Y_1(t)=0\quad{\rm on}\;\; (0,1].
\end{equation}
Since $u_\infty$ is bounded in $(0,1]\ti B_2(\bar x)$ by (\ref{C6}),  $X_1(t)$ and $Y_1(t)$ remains bounded on $(0,1]$.
$$Z_1(t)=\myint{B_1}{}\!\!\!u_\infty(t,x)(-\Gd)^s\eta(x) dx+\myint{B^c_1}{}\!\!\!u_\infty(t,x)(-\Gd)^s\eta(x) dx:=Z_{1,1}(t)+Z_{1,2}(t).
$$
Since $\eta$ has its support in $B_1(\bar x)$, there exists $c_8>0$ such that
$$(-\Gd)^s\eta(x)\leq -c_8\qquad\forall\, x\in B_1(0).
$$
Using (\ref{3.6}) we obtain that
$$\myint{0}{1}Z_{1,1}(s)ds=-\infty.
$$
Using (\ref{II-15}) in \rlemma{test}, we have
$$\BA {lll}\abs{Z_{1,2}(t)}\leq c_5\norm\eta_{C^2}\myint{B_1^c}{}\myfrac{u_\infty(t,x) dx}{1+\abs{x}^{N+2s}}\leq c_9\myint{1}{\infty}\myfrac{r^{N-1}dr}{1+r^{N+2s}}.\EA$$
Hence
$$\myint{0}{1}Z_{1}(s)ds=-\infty.
$$
Integrating (\ref{C7}) it contradicts the boundedness of $X_1$ and $Y_1$. Hence, for any $z\in\BBR^N$,
\begin{equation}\label{C8}
\displaystyle\limsup_{t\to 0}u_\infty(t,z)=\infty.
\end{equation}
Using again the fact that $x\mapsto u_\infty(t,x)$ is radial and decreasing with respect to $\abs x$, we get
\begin{equation}\label{C9}
\displaystyle\limsup_{t\to 0}\myint{B_\gr(z)}{}u_\infty(t,x)dx=\infty\qquad\forall\,\gr>0.
\end{equation}
By Theorem C, we infer that $u_\infty(t,x)\geq u_{z,\infty}(t,x)=u_{\infty}(t,x-z)$. Interchanging $0$ and $z$ we conclude again that $u_\infty(t,x)$ depends only on $t$, hence it coincides with $U_{\gb,p}(t)$, and clearly $\CS_{u_\infty}=\BBR^N$.\qeda

\subsection{Proof of Theorem G  (ii)}
We assume that $\gk\geq 1$ and $\mathcal{L}=\{0_{\BBR^\gk}\}\times\R^{N-\gk}$. We set $x=(x_1,x')\in \BBR^\gk\ti\BBR^{N-\gk}$.
We use Theorem G  $(i)$ with $N$ replace by $N-\gk$ to prove the part $(ii)$.  If $x=(x_1,x')$, then $\bar x=(0,x')\in\CL$, hence by \cite[Th 1.3 (ii)]{CVW}
$$\BA {lll}
u_\infty(t,x-\bar x)\geq \myfrac{c_{10}t^{-\frac{N+2s}{2s}}}{1+(t^{-\frac{1}{2s}}\abs{x-\bar x})^{N+2s}}=
\myfrac{c_{10}t^{-\frac{N+2s}{2s}}}{1+(t^{-\frac{1}{2s}}\abs{x_1})^{N+2s}}.
\EA
$$
By Theorem C, we obtain
\begin{equation}\label{C10}
u(t,x)\geq \myfrac{c_{10}t^{-\frac{N+2s}{2s}}}{1+(t^{-\frac{1}{2s}}\abs{x_1})^{N+2s}}\qquad\forall\, (t,x):=(t,x_1,x')\in \BBR_+^*\ti  \BBR^\gk\ti\BBR^{N-\gk}.
\end{equation}
For $n\in\BBN^*$, set $t_n=n^{-2s}$, $\gr_0=\int_{\BBR^{N-\gk}}\frac{dy'}{1+\abs{y'}^{N+2s}}$, $\gl_0=\gr_0^{\frac{1}{p-1}}$ and
$$f_n(x_1)=\myfrac{c_{10}\gl_0(\gr_0t_n)^{-\frac{1+\gb}{p-1}}}{1+\left((\gr_0t_n)^{-\frac{1}{2s}}\abs{x_1}\right)^{N+2s}}.
$$
Then
$$\myint{\BBR^\gk}{}f_n(x_1)dx_1=c_\gk t_n^{\frac{\gk}{2s}-\frac{1+\gb}{p-1}}\myint{0}{\infty}\myfrac{r^{\gk-1}dr}{1+r^{N+2s}}
$$
 for some $c_\gk>0$. Since $p<1+\frac{2s(1+\gb)}{\gk}$ the above integral is finite for any $n$ but tends to $\infty$ with $n$. Hence we fix $m>0$, then for any $n\in\BBN^*$ there exists $\ge_{n,m}>0$ such that
$$\myint{\abs {x_1}\leq \ge_{n,m}}{}f_n(x_1)dx_1=m=c_\gk t_n^{\frac{\gk}{2s}-\frac{1+\gb}{p-1}}\myint{0}{\ge_{n,m}t_n^{-\frac{N+2s}{2s}}}\myfrac{r^{\gk-1}dr}{1+r^{N+2s}}.
$$
Hence $\ge_{n,m}\to 0$ when $n\to\infty$. This implies that for any $\gz\in C^\infty_0(\BBR^\gk)$,
$$\displaystyle\lim_{n\to\infty}\myint{\abs {x_1}\leq \ge_{n,m}}{}f_n(x_1)dx_1= m\gz(0).
$$
Equivalently $f_{n,m}:=f_n\chi_{_{B_{\ge_{n,m}}}}\to m\gd_0$ in the sense of measures in $\BBR^\gk$.
Let $w_{n,m}$ be the solution of
\begin{equation}\label{C11}\BA {lll}
\prt _tu+(-\Gd)_{_{\BBR^\gk}}^s u+t^\gb u^p=0\qquad&{\rm in}\;\; \BBR^*_+\ti\BBR^\gk\\[2mm]
\phantom{dd-----\ \ \ }
\!u(0,\cdot)=f_{n,m}&{\rm in}\;\; \BBR^\gk,
\EA\end{equation}
in which formula $(-\Gd)_{_{\BBR^\gk}}^s$ denotes the fractional Laplacian in $\BBR^\gk$, an index omitted if $\gk=N$. Then
$\tilde w_{n,m}(t,x_1,x')=w_{n,m}(t,x_1)$ is a solution of
\begin{equation}\label{C12}\BA {lll}
\prt _tu+(-\Gd)^s u+t^\gb u^p=0\qquad&{\rm in}\;\; \BBR^*_+\ti\BBR^N\\[1mm]
\phantom{ddd^,-\ \ -,,-u}
u(0,\cdot)=\tilde f_{n,m}&{\rm in}\;\; \BBR^N,
\EA\end{equation}
with $\tilde f_{n,m}(x_1,x')=f_{n,m}(x_1)$. Since $u(t_n,x)\geq \tilde f_n(x)$ in $\BBR^N$, we obtain by  the comparison principle that
$u(t+t_n,x)\geq \tilde w_{n,m}(t,x)$ in $\BBR^*_+\ti \BBR^N$. Hence, by letting successively $n\to\infty$ and $m\to\infty$,
\begin{equation}\label{C13}
u(t+t_n,x_1,x')\geq w_{n,m}(t,x_1)\Longrightarrow u(t,x_1,x')\geq u^\gk_{m\gd_0}(t,x_1)\Longrightarrow u(t,x_1,x')\geq u^\gk_{\infty}(t,x_1),
\end{equation}
where we have denoted by $u^\gk_{m\gd_0}$ and $u^\gk_{\infty}$ respectively the solution of the equation in (\ref{C11}) with $m\gd_0$ as initial data and the limit of this solution when $m\to\infty$.  Since $1<p\leq 1+\frac{2s}{\gk+2s}$, $u^\gk_{\infty}=U_{p,\gb}$ by (i), which ends the proof.\qeda\medskip

\nind\Remark It appears interesting to investigate whether the fact that the singular set $\CS_u$ contains a $(N$-$\gk)$-dimensional plane can be replaced by $\CS_u$ contains a $(N$-$\gk)$-dimensional submanifold.


\section{Solution with a given initial trace: the general case }
\setcounter{equation}{0}

\subsection{Problem with initial data measure}

If $\gn$ is a bounded Radon measure on an open set $\CR\subset\BBR^N$, that we note $\gn\in\mathfrak M^b(\CR)$), we extend it by $0$ in $\CS=\BBR^N\setminus \CR$ and the new measure still denoted by $\gn$, belongs to the space $\mathfrak M^b(\BBR^N)$ of bounded Radon measures on $\BBR^N$. Conversely, if $\gn\in \mathfrak M^b(\BBR^N)$ vanishes  on $\CS=\BBR^N\setminus \CR$ its restriction to $\CR$ belongs to $\mathfrak M^b(\CR)$.

\begin{definition}\label{good} A nonnegative bounded Radon measure $\gn$ in $\BBR^N$ is an admissible measure if there exists a function $u=u_\gn\in L^1(Q_T)$ with $ t^\gb u^p\in L^1(Q_T)$  solution of 
\begin{equation}\label{g-1}
\begin{array}{lll}
\partial_t u + (-\Delta)^s  u+ t^\beta u^p=0\quad & {\rm in}\quad Q_\infty\\[1mm]
\phantom{-------}
u(0,.)=\nu.
\end{array}
\end{equation}
It is a good measure if the sequence $u_{\gn_n}$ of solutions of (\ref{g-1}) with initial data $\gn_n=\gn\ast\gr_n$, where $\{\gr_n\}$ is a sequence of mollifiers, converges to $u_{\gn}$ in $L^1(Q_T)$ and if $ t^\gb u_{\gn_n}^p$ converges to $ t^\gb u^p\in L^1(Q_T)$.
\end{definition}
Uniqueness of solutions is proved in \cite{CVW} as a result of the choice of $\BBY_{s,T }$ as space of test functions. Notice also that if $p<p^*_\gb$ any nonnegative bounded measure is good. The following result will be useful in the sequel.

\bprop{add} Let $p>1$ and $\gb>-1$. If $\gn,\gm\in \mathfrak M_+^b(\BBR^N)$ are good measures (resp. admissible measures), then $\gn+\gm$ is a good measure (resp. admissible measure).
\es
\Proof We set $\gn_n=\gn\ast\gr_n$ and $\gm_n=\gm\ast\gr_n$ and denote by  $u_{\gn_n}$, $u_{\gm_n}$ and $u_{\gn_n+\gm_n}$ the solutions of the initial value problem (\ref{g-1}) with $\gn$ replaced by $\gn_n$, $\gm_n$ and $\gn_n+\gm_n$ respectively. Since $p>1$,  $u_{\gn_n}+u_{\gm_n}$ is a supersolution of (\ref{In-15}). Hence $u_{\gn_n+\gm_n}\leq u_{\gn_n}+u_{\gm_n}$. When $n\to\infty$, $u_{\gn_n+\gm_n}$ converges a.e. to some function $u$ (see \cite{CVW}). Since $u_{\gn_n}$ and $u_{\gm_n}$ converges in $L^1(Q_T)$, the sequence $u_{\gn_n+\gm_n}$ is uniformly integrable in $Q_T$, it converges to some $w$ (up to extraction of a subsequence). Furthermore,
$$(u_{\gn_n+\gm_n})^p\leq (u_{\gn_n}+u_{\gm_n})^p\leq 2^{p-1}\left(u^p_{\gn_n}+u^p_{\gm_n}\right).$$
Since $t^\gb u^p_{\gn_n}$ and $t^\gb u^p_{\gm_n}$ converges in $L^1(Q_T)$ to $t^\gb u^p_{\gn}$ and $t^\gb u^p_{\gm}$ respectively, they are uniformly integrable. Hence the sequence $\{t^\gb (u_{\gn_n+\gm_n})^p\}$ is uniformly integrable in $Q_T$ and thus, up to extraction of a second subsequence, $t^\gb (u_{\gn_n+\gm_n})^p$ converges to $t^\gb w^p$ in $L^1(Q_T)$. Going to the limit in the  formulation $(\ref{II-16+})$ of the fact $u_{\gn_n+\gm_n}$ is a weak solution of (\ref{g-1}) with initial data $\gn_n+\gm_n$, it follows that $w$ satisfies the same equation (\ref{g-1}) but now with initial data $\gm+\gn$.
By uniqueness (see \cite[Th 1.1]{CVW} and notice that therein uniqueness  needs no more condition on $h$ than monotonicity), $w=u_{\gn+\gm}$ and the whole sequence $\{u_{\gn_n+\gm_n}\}$ converges to $u_{\gn+\gm}$. The proof in the case of admissible measures is similar. \qeda\medskip

\bprop{sup} Let $p>1$ and $\gb>-1$. If $\{\gn_k\}\subset \mathfrak M_+^b(\BBR^N)$ is a nondecreasing sequence of admissible measures converging to $\gn\in\mathfrak M_+^b(\BBR^N)$, then $\gn$ is an admissible measure.
\es
\Proof The sequence $\{u_{\gn_k}\}$ is increasing. Furthermore,
$$u_{\gn_k}\leq \BBH_s[\gn_k]\leq \BBH_s[\gn].
$$
Hence there exists some $u\in L^1(Q_T)$ for any $T>0$, $u\geq 0$, such that $u_{\gn_n}\to u$ in  $L^1(Q_T)$ and a.e. in $Q_\infty$.  By identity (3.25) in the proof of \cite[Th. 1.1]{CVW}, we have for $\gt\geq T$,
\bel{h-2-}
\myint{0}{1}\myint{\BBR^N}{}\left(u_{\gn_k}+(\gt-t)t^\gb u^q_{\gn_k}\right)dxdt+(\gt-T)\myint{\BBR^N}{}u_{\gn_k}(T,x)dx=\gt\myint{\BBR^N}{}d\gn_k
\leq 2\gt\norm\gn_{\mathfrak M_+^b}.
\ee
Hence $t^\gb u^q\in L^1(Q_T)$ and $t^\gb u^q_{\gn_k}\to t^\gb u^q$ in $L^1(Q_T)$. By (\ref{II-16+}) there holds
$$\myint{0}{1}\myint{\BBR^N}{}\left[u_{\gn_k}\left(-\prt_t\xi+(-\Gd)^s\xi\right)+t^\gb u^q_{\gn_k}\xi\right]dxdt=
\myint{\BBR^N}{}\xi(0,x)d\gn_k(x),
$$
for all $\xi\in \BBY_{s,T }$, so it follows that $u=u_\gn$. Hence $\gn$ is an admissible measure. \qeda

\medskip

The whole description of the set of admissible measures necessitates the introduction of Bessel capacities as in the case $s=1$, see \cite{MV1}, \cite{MV5}. We have a first partial answer.

\blemma{add-2} Let $p>1$ and $\gb>-1$. If $\gn\in \mathfrak M_+^b(\BBR^N)$ satisfies  $t^\gb(\BBH_s[\gn])^p\in L^1(Q_T)$, then $\gn$ is a good measure.
\es
\Proof Let $\gn_n=\gn\ast\gr_n$, by the maximum principle
$$u_{\gn_n}\leq \BBH_s[\gn_n]=\BBH_s[\gn\ast\gr_n]=\BBH_s[\gn]\ast\gr_n.$$
Since $\BBH_s[\gn]\in L^1(Q_T)$, $\BBH_s[\gn]\ast\gr_n\to \BBH_s[\gn]$ in $L^1(Q_T)$. Similarly $t^\gb(\BBH_s[\gn]\ast\gr_n)^p\to t^\gb(\BBH_s[\gn])^p$ in $L^1(Q_T)$. Since $u_{\gn_n}\leq \BBH_s[\gn]\ast\gr_n$, we conclude that the sequences
$\{u_{\gn_n}\}$ and $\{t^\gb u^p_{\gn_n}\}$ are uniformly integrable in $L^1(Q_T)$, hence they are precompact by Vitali's convergence theorem. As in the proof of \rprop{add} any cluster point $w$ in the $L^1(Q_T)$-topology of the sequence $\{u_{\gn_n}\}$
 is a weak solution of (\ref{g-1}) with initial data $\gn$ hence $w=u_\gn$ and $u_{\gn_n}\to u_\gn$ by uniqueness of the solution. \qeda

\bprop{add-3} Let $p>1$ and $\gb>-1$. Assume that $\gn\in \mathfrak M_+^b(\BBR^N)$ with Lebesgue decomposition $\gn=\gn_0+\gn_s$, where
$\gn_0$ and $\gn_s$ belong to $\mathfrak M_+^b(\BBR^N)$, $\gn_0\in L^1(\BBR^N)$ and $\gn_s$ is singular with respect to the $N$-dim Lebesgue measure. If $t^\gb(\BBH_s[\gn_s])^p\in L^1(Q_T)$, then $\gn$ is a good measure.
\es
\Proof By \cite[Lemma 3.2]{CVW} there exists a unique solution $u_{\gn_0}$ (resp. $\gn_s$) to problem (\ref{g-1}) with $\gn$ replaced by $\gn_0$.
(resp. $\gn_s$). By \cite[Lemma 3.2]{CVW} the sequences $\{u_{\gn_0\ast\gr_n}\}$ and $\{t^\gb u^p_{\gn_0\ast\gr_n}\}$ are Cauchy sequences in $L^1(Q_T)$. By
\rlemma{add-2}, the sequences $\{u_{\gn_s\ast\gr_n}\}$ and $\{t^\gb u^p_{\gn_s\ast\gr_n}\}$ share the same property. Hence $\gn_0$
and $\gn_s$ are good measures and we conclude with \rprop{add}. $\phantom{------}$\qeda\medskip

We recall some classical results about Bessel potentials, capacities and interpolation.
For $0<\gg<N$, the Bessel kernel $J_\gg$ is defined in $\BBR^N\setminus\{0\}$ by $J_\gg(x)=\CF^{-1}\left((1+|\xi|^2)^{-\frac{\gg}{2}}\right)$, where $\CF$ is the Fourier transform in $\BBR^N$, and the Bessel potential of a positive measure is
\bel{h-2}{\bf J}_\gg[\gm](x)=J_\gg\ast\gm(x)=\myint{\BBR^N}{}J_\gg(x-y)d\gm(y).
\ee
For $1\leq r<\infty$, the Bessel capacity ${\rm cap\,}^{\BBR^N}_{\gg,r}$  of a compact set is
\bel{h-3}{\rm cap\,}^{\BBR^N}_{\gg,r}(K)=\inf\{\norm{{\bf J}_\gg[\gz]}^r_{L^r(\BBR^N)}:\gz\in\gw_K\},
\ee
where $\gw_K$ is the subset of nonnegative  function belonging to the Schwartz space $\CS(\BBR^N)$ , with value larger or equal to $1$ on $K$. Furthermore
\bel{h-4}(-\Gd+I)^{\frac\gg 2}\gf=\gz\Longleftrightarrow {\bf J}_\gg[\gz]=\gf.\ee

If a linear m-accretive operator $A$ in $L^r(\BBR^N)$ with domain $D(A)$ is the infinitesimal generator of an analytic semigroup of bounded linear operators
$S^A(t)$, i.e.
\bel{h-5}u(t)=S_t^{-A}\gu\quad\forall t\geq 0\Longleftrightarrow\frac{du}{dt}+Au=0\quad{\rm on}\;\BBR_+^*\,,\;u(0)=\gu,
\ee
the real interpolation classes between $D(A)$ and $L^r(\BBR^N)$ can be obtained (see \cite[p. 96]{Tri}) by
\bel{h-1}\left[D(A),L^r(\BBR^N)\right]_{\gth,r}=\left\{\gu\in L^r(\BBR^N):\myint{0}{1}\myint{\BBR^N}{}t^{(1-\gth)r}|AS_t^A\gu|^{r}dx\myfrac{dt}{t}<+\infty\right\},
\ee
and
\bel{h-6}\norm\gu_{\left[D(A),L^r\right]_{\gth,r}}\equiv \norm\gu_{L^r}+\left(\myint{0}{1}\norm{t^{1-\gth}AS_t^A\gu}^r_{L^r}\myfrac{dt}{t}\right)^{\frac 1r}.
\ee
If $A=(-\Gd)^s+ I$, its domain $D(A)$ in $L^r(\BBR^N)$ is the Bessel potential space $(I-\Gd)^{-s}(L^r(\BBR^N))=L^{2s,r }(\BBR^N)$: the result is stated in \cite[Th 1]{Gur} but it is an easy consequence of \cite[Chap. 4, Th 3]{Stein} applied to the Fourier multipliers $(I+\abs {\xi}^2)^s(1+\abs {\xi}^{2s})^{-1}$ and $(I+\abs {\xi}^2)^{-s}(1+\abs {\xi}^{2s})$. By classical interpolation properties of Bessel potential spaces (see e.g. \cite{Tri}),

\bel{h-7}\left[D(A),L^r(\BBR^N)\right]_{\gth,r}=L^{2\gth s,r}(\BBR^N)=(I-\Gd)^{-\gth s}(L^r(\BBR^N)).
\ee
Since $A$ is coercive, \cite[Sec. 1.14.5]{Tri}, (\ref{h-6}) can be replaced by

\bel{h-6+}\norm\gu_{\left[D(A),L^r\right]_{\gth,r}}\equiv \left(\myint{0}{1}\norm{t^{1-\gth}AS_t^A\gu}^r_{L^r}\myfrac{dt}{t}\right)^{\frac 1r}.
\ee

\bprop{adm} Let $N\geq 1$, $p>1$ and $-1<\gb<p-1$. If problem (\ref{g-1}) admits a positive solution $u_\gn$ for some $\gn\in \mathfrak M_b^+(\BBR^N)$, then $\gn$ vanishes on Borel subsets of $\BBR^N$ with zero Bessel capacity ${\rm cap\,}^{\BBR^N}_{\frac{2s(1+\gb)}{p},p'}$ where $p'=\frac{p}{p-1}$, i.e.
\begin{equation}\label{g-3}
\forall K\subset\BBR^N,K\;\;{\rm Borel},\,{\rm cap\,}^{\BBR^N}_{\frac{2s(1+\gb)}{p},p'}(K)=0\Longrightarrow \gn(K)=0.
\end{equation}
\es
\Proof Assume that $u:=u_\gn$ is the solution of (\ref{g-1}). Since ${\rm cap\,}^{\BBR^N}_{\frac{2s(1+\gb)}{p},p'}$ is a Choquet capacity, let $K\subset\BBR^N$ is compact and   $\gz\in \CS(\BBR^N)$ be such that $0\leq\gz$ in $\BBR^N$ and $\gz\geq  1$ on $K$. We set $\Phi=e^{- t}\BBH_s[\gz]$ and take $\Phi^{p'}$ as a test function. Then
\bel{h8}\myint{\BBR^N}{}u\Phi^{p'}(1,.)dx+\myint{0}{1}\myint{\BBR^N}{}\left[\left(-\prt_t\Phi^{p'}+(-\Gd)^s\Phi^{p'}\right)u+t^\gb u^p\Phi^{p'}\right] dxdt=\myint{\BBR^N}{}\Phi^{p'}d\gn.
\ee
Note that $(-\Gd)^s\Phi^{p'}\geq p'\Phi^{p'-1}(-\Gd)^s\Phi$ and $\prt_t\Phi+(-\Gd)^s\Phi+\Phi=0$, hence
$$-\prt_t\Phi^{p'}+(-\Gd)^s\Phi^{p'}\geq 2p'\Phi^{p'-1}(-\Gd)^s\Phi.
$$
Then we adapt the duality argument of \cite{BP} and get from H\"older's inequality
$$\BA {lll}\myint{\BBR^N}{}u\Phi^{p'}(1,x)dx+2p'\left(\myint{0}{1}\myint{\BBR^N}{}t^{-\frac{p'\gb}{p}}\left|(-\Gd)^s\Phi+\Phi\right|^{p'}dx dt\right)^{\frac{1}{p'}}
\left(\myint{0}{1}\myint{\BBR^N}{}
t^\gb u^p\Phi^{p'}dxdt\right)^{\frac{1}{p}}\\[4mm]
\phantom{----))------------------}
+\myint{0}{1}\myint{\BBR^N}{}
t^\gb u^p\Phi^{p'}dxdt\geq \gn(K).
\EA$$
Applying (\ref{h-1}), (\ref{h-7}) with $r=p'$, $\gth=\frac{1+\gb}{p}$, we obtain directly for some $c_{11}>1$,
\bel{h9}\BA {lll}\myfrac{1}{c_{11}}\norm\gz_{L^{\frac{2s(1+\gb)}{p},p'}}\leq
\left(\myint{0}{1}\myint{\BBR^N}{}t^{-\frac{p'\gb}{p}}\left|(-\Gd)^s\Phi+\Phi\right|^{p'}dx dt\right)^{\frac{1}{p'}}
\leq c_{11}\norm\gz_{L^{\frac{2s(1+\gb)}{p},p'}}.
\EA\ee
If ${\rm cap\,}^{\BBR^N}_{\frac{2s(1+\gb)}{p},p'}(K)=0$, there exists a sequence $\{\gz_n\}\subset \gw_K$ such that $\norm{\gz_n}_{L^{\frac{2s(1+\gb)}{p},p'}}\to 0$ as $n\to\infty$. Furthermore it is possible to assume $\gz_n\leq 1$  in $\BBR^N$ (see \cite{AdPo}). Hence, up to a subsequence, $\gz_n\to 0$ a.e. in $\BBR^N$. This implies $\Phi_n\leq 1$ and $\Phi_n\to 0$ a.e. in $Q_\infty$. Therefore,
$$\displaystyle\lim_{n\to\infty}\myint{\BBR^N}{}u\Phi_n^{p'}(1,x)dx
\;\;{\rm and}\quad \lim_{n\to\infty}\myint{0}{1}\myint{\BBR^N}{}
t^\gb u^p\Phi_n^{p'}dxdt=0.
$$
Combining the previous inequalities we infer that $\gn(K)=0$.\qeda\medskip

\rprop{adm} is the necessary condition in Theorem H. The next result provides the sufficient condition.

\bprop{adm2}Let $N\geq 1$, $p>1$, $-1<\gb<p-1$ and $\gn\in \mathfrak M_b^+(\BBR^N)$ which vanishes on Borel subsets of $\BBR^N$ with zero  ${\rm cap}^{\BBR^N}_{\frac{2s(1+\gb)}{p},p'}$-Bessel  capacity. Then $\gn$ is an admissible measure.
\es
\Proof If $\gn$ vanishes Borel subsets with zero  ${\rm cap}^{\BBR^N}_{\frac{2s(1+\gb)}{p},p'}$, there exists an increasing sequence of nonnegative measures
$\{\gn_n\}\subset \left(L^{{\frac{2s(1+\gb)}{p-1},p'}}(\BBR^N)\right)'=L^{-{\frac{2s(1+\gb)}{p-1},p}}(\BBR^N)$ such that $\gn_n\to\gn$ in the sense of measures. This results is classical in the integer case and a proof in the Bessel case (similar in fact) can be found in \cite[Prop. 3.6]{Ve1}. \\
Next let $\gz\in \CS(\BBR^N)$ and $\Gf=e^{- t}\BBH_s[\gz]$, then
$$\myint{\BBR^N}{}\Gf\BBH_s[\gn_n](1,x)dx+\myint{0}{1}\myint{\BBR^N}{}\BBH_s[\gn_n]\left(2(-\Gd)^s\Gf +\Gf \right)dxdt=\myint{\BBR^N}{}\gz d\gn_n.
$$
Hence
$$\myint{0}{1}\myint{\BBR^N}{}\BBH_s[\gn_n]\left((-\Gd)^s\Gf +\Gf\right)dxdt\leq \norm{\gn_n}_{L^{-{\frac{2s(1+\gb)}{p-1},p}}}\norm\gz_{L^{{\frac{2s(1+\gb)}{p-1},p'}}}.
$$
Consider the mapping
$$\gz\mapsto L(\gz)=\myint{0}{1}\myint{\BBR^N}{}t^{\frac\gb p}\BBH_s[\gn_n]t^{-\frac\gb p}\left((-\Gd)^s\Gf +\gd\Gf\right)dxdt.
$$
It satisfies
\bel{h10}\BA{lll}\abs{L(\gz)}\leq \norm{\gn_n}_{L^{-{\frac{2s(1+\gb)}{p-1},p}}}\norm\gz_{L^{{\frac{2s(1+\gb)}{p-1},p'}}}\\[4mm]
\phantom{\abs{L(\gz)}}
\leq
c\norm{\gn_n}_{L^{-{\frac{2s(1+\gb)}{p-1},p}}}\left(\myint{0}{1}\myint{\BBR^N}{}t^{-\frac{p'\gb}{p}}\left|(-\Gd)^s\Phi+\gd\Phi\right|^{p'}dx dt\right)^{\frac{1}{p'}},
\EA\ee
by (\ref{h9}). Hence $t^{\frac\gb p}\BBH_s[\gn_n]\in L^p(Q_1)$ and
\bel{h11}
\left(\myint{0}{1}\myint{\BBR^N}{}t^{\gb}\left(\BBH_s[\gn_n]\right)^pdx dt\right)^{\frac{1}{p}}\leq
c_{12}\norm{\gn_n}_{L^{-{\frac{2s(1+\gb)}{p-1},p}}}.
\ee
Hence the $\gn_n$ are good measures by \rlemma{add-2}. Then by  \rprop{sup}, $\gn$ is an admissible measure.
$\phantom{----------}$\qeda\medskip

\nind\Remark  When $s=1$ and $\gb=0$, it is proved in \cite{BP} that the admissibility condition for measures is strongly linked to the removability for Borel sets in the sense that  if $K\subset\BBR^N$ is a Borel set with zero ${\rm cap}^{\BBR^N}_{\frac{2}{p},p'}$-capacity,  any $u\in C(\overline Q_\infty\setminus\{(0,K)\})\cap C^{1,2}(Q_\infty)$ solution of (\ref{In-10}) in $Q_\infty$ which vanishes on $(0,x)$ for any $x\in\BBR^N\setminus K$ is identically zero. The set $K$ is said {\it removable}. Furthermore, the condition is also necessary. Now, for equation (\ref{In-10}) it is clear that a compact set $K$ with positive ${\rm cap}^{\BBR^N}_{\frac{2s(1+\gb)}{p},p'}$-capacity it is not removable since it is the support  of the capacitary measure (a positive measure belonging to $L^{\frac{2s(1+\gb)}{p},p'}(\BBR^N)$, \cite[Chap. 2]{AdHe}) which is a good measure by \rlemma{add-2}. {\it We conjecture that the condition ${\rm cap}^{\BBR^N}_{\frac{2s(1+\gb)}{p},p'}(K)=0$ implies the removability of the compact set $K$ for equation (\ref{In-15}) in the sense given above}. 

\subsection{Barrier function for $N=1$}

We set
\bel{bn-1}W(z)=\left\{
\begin{array}{lll}
\myfrac{\ln(e+z^2)}{1+z^{1+2s}}\quad & {\rm if}\quad z\geq  0\\[4mm]
1\quad & {\rm if}\quad z<  0,
\end{array}
\right.
\ee
where  $e$ is Neper constant, and
\begin{equation}\label{17-09-0}
 w(t,x)= t^{-\frac{1+\beta}{p-1}} W(t^{-\frac1{2s}}x),\quad \forall\, (t,x)\in \BBR^*_+\times\R.
 \end{equation}
When $t\to 0$, the function $w$ satisfies
 \bel{bn-2}\BA {lll}
 (i)\qquad &w(t,x)= \myfrac{2t^{\frac{1+2s}{2s}-\frac{1+\beta}{p-1}}\ln t}{|x|^{1+2s}} (1+o(1))&\qquad{\rm if }\quad x>0,\\[3mm]
 (ii)\qquad &w(t,x)= t^{-\frac{1+\beta}{p-1}}&\qquad{\rm if }\quad x\leq 0.
 \EA\ee

\begin{lemma}\label{lm 14-09-0}
Assume that $p>1+\frac{2s(1+\beta)}{1+2s}$. Then there exists $\lambda_0>0$ such that for $\lambda\geq \lambda_0$, the function
$w_\gl:=\gl w$ satisfies
\bel{14-09-0}\BA {lll}\displaystyle
\partial_tw_\gl+(-\Gd)_{_{\BBR}}^s w_\gl+t^\beta w_\gl^p\geq 0\quad&{\rm in}\quad   \BBR^*_+\times\R\\[2mm]
\phantom{------,,}
\displaystyle\lim_{t\to 0}w(t,x)=0\qquad&{\rm if }\quad x>0\\[2mm]
\phantom{------,,}
\displaystyle\lim_{t\to 0}w(t,x)=\infty\qquad&{\rm if }\quad x\leq 0.
 \EA
\ee
\end{lemma}

\nind\Proof
Clearly the assertions concerning the limit of $w(t,x)$ when $t\to 0$ are satisfied since  $\frac{1+2s}{2s}-\frac{1+\beta}{p-1}>0$ by assumption.
Then
$$
\partial_tw_\lambda(t,x) = -\frac{\lambda(1+\beta)}{p-1}t^{-\frac{1+\beta}{p-1}-1}w(z)
-\frac{\lambda}{2s}t^{-\frac{1+\beta}{p-1}-1}w'(z)z,
$$
with $z=t^{-\frac{1}{2s}}x$ and
$$(-\Delta)_1^s w_\lambda(t,x)=\lambda t^{-\frac{1+\beta}{p-1}-1}(-\Delta)_1^s w(z).$$
Hence
\begin{equation} \label{15-09-00}
\arraycolsep=1pt
\begin{array}{lll}
\partial_t w_\lambda(t,x)+(-\Gd)_{_{\BBR}}^s w_\lambda(t,x)+t^\beta w^p_\lambda(t,x)\\[2mm]
\phantom{----}
   =\lambda t^{-\frac{1+\beta}{p-1}-1}\left[(-\Gd)_{_{\BBR}}^s w(z)-\myfrac{1}{2s}w'(z)z-\myfrac{1+\beta}{p-1}w(z)+\lambda^{p-1}w^p(z)\right].
   \end{array}
\end{equation}
If $z>0$, we obtain that
$$ \arraycolsep=1pt
\begin{array}{lll}
-\myfrac{1}{2s}w'(z)z-\myfrac{1+\beta}{p-1}w(z) = \left[\myfrac{1+2s}{2s}\myfrac{z^{1+2s}}{1+z^{1+2s}}-\myfrac{1+\beta}{p-1}
-\myfrac{z^2(e+z^2)^{-1}}{s\ln(e+z^2)}\right]w(z).
\end{array}
$$
Since $\frac{1+2s}{2s}>\frac{1+\beta}{p-1}$, $\displaystyle\lim_{z\to\infty}$$\frac{z^{1+2s}}{1+z^{1+2s}}=1$ and $\displaystyle\lim_{z\to\infty}$$\frac{1}{\ln(e+z^2)}=0$, then
there exist $R_0>0$ and $\sigma_0>0$ such that
\begin{equation}\label{14-09-5}
-\frac1{2s}w'(z)z-\frac{1+\beta}{p-1}w(z)\geq  \sigma_0 w(z)\qquad \forall\, z\geq  R_0.
\end{equation}
Next we deal with $(-\Gd)_{_{\BBR}}^s w(z)$ and put

$$\tilde w(z)=\frac{\ln(e+z^2)}{1+ |z|^{1+2s}}\qquad \forall\, z\in\R, $$
so that $(-\Gd)_{_{\BBR}}^s w(z)=(-\Gd)_{_{\BBR}}^s \tilde w(z)+(-\Delta)_1^s (1-\tilde w\chi_{_{\R_-}})(z)$.\\
 For $z>2$, using the equivalent definition of fractional Laplacian, we have that
\begin{equation} \label{15-09-0}
\arraycolsep=1pt
\begin{array}{lll}
-(-\Gd)_{_{\BBR}}^s \tilde w(z)=\myfrac{a_{1,s}}{2}\myint{-\infty}{\infty}\myfrac{\frac{\ln(e+|z+\tilde y|^2)}{1+|z+\tilde y|^{1+2s}}
+\frac{\ln(e+|z-\tilde y|^2)}{1+|z-\tilde y|^{1+2s}}-\frac{2\ln(e+z^2)}{1+z^{1+2s}}}{|\tilde y|^{1+2s}}d \tilde y\\[4mm]
\phantom{------\ }
   =\myfrac{a_{1,s}w(z)}{2z^{2s}}\myint{-\infty}{\infty}\myfrac{I_z(y)}{|y|^{1+2s}}dy,
   \end{array}
\end{equation}
where
$$\arraycolsep=1pt
\begin{array}{lll}
 I_z(y) = \myfrac{1+z^{N+2s}}{1+z^{1+2s}|1+y|^{1+2s}}\myfrac{\ln(e+z^2|1+y|^2)}{\ln(e+z^2)} 
  +\myfrac{1+z^{1+2s}}{1+z^{1+2s}|1-y|^{1+2s}} \myfrac{\ln(e+z^2|1-y|^2)}{\ln(e+z^2)}-2.
 \end{array}
$$

\nind {\it Step 1: There exists $c_{12}>0$ such that}
 \begin{equation}\label{15-09-1}
   \myint{\frac12\leq  |y|\leq \frac32}{}\myfrac{I_z(y)}{|y|^{1+2s}}dy\leq  \frac{ c_{12}}{w(z)z}.
 \end{equation}
Actually, for $-\frac32 <y<-\frac12$, there exists $c_{13}>0$ such that
$$\frac{1+z^{1+2s}}{1+z^{1+2s}|1-y|^{1+2s}} \frac{\ln(e+z^2|1-y|^2)}{\ln(e+z^2)}\leq  c_{13}$$
and  then
\begin{eqnarray*}
 \int_{-\frac32}^{-\frac12}\frac{I_z(y)}{|y|^{1+2s}}dy &\leq & 2\int_0^{\frac12}\frac{1+z^{1+2s}}{1+(z r)^{1+2s}}\frac{\ln(e+z^2 r^2)}{\ln(e+z^2)}dr+c_{14}
    \\&\leq &\frac{2}{w(z)z}\int_0^{\infty}\frac{\ln(e+t^2)}{1+ t^{1+2s}}dt+c_{14}
    \\&\leq &\frac{c_{15}}{w(z)z},
\end{eqnarray*}
where $c_{14},c_{15}>0$, and the last inequality holds since $w(z)z\to0$ as $z\to+\infty$.
Similarly,
$$\int_{\frac12}^{\frac32}\frac{I_z(y)}{y^{N+2s}}dy_1\leq  \frac{c_{16}}{w(z)z}.$$

\nind {\it Step 2: There exists $c_{17}>0$ such that}
 \begin{equation}\label{15-09-2}
   \int_{-\frac12}^{\frac12}\frac{I_z(y)}{|y|^{1+2s}}dy\leq  c_{17}.
 \end{equation}
Indeed, since function $I_z$ is $C^2$ in $[-\frac12,\frac12]$ and satisfies
$$I_z(0)=0\;\;{\rm and}\;\;  I_z(y)=I_z(-y),$$
then $I_z'(0)=0$ and there exists $c_{18}>0$ such that
$$ |I_z''(y)|\leq  c_{18} \quad{\rm for\ any\ } y\in [-\textstyle\frac12,\textstyle\frac12]. $$
Then we have
$$|I_z(y)|\leq  \frac{c_{18}}2y^2 \quad{\rm for\ any\ } y\in[-\textstyle\frac12,\textstyle\frac12],$$
which implies that
$$
\left|\int_{-\frac12}^{\frac12}\frac{I_z(y)}{|y|^{1+2s}}dy\right|\leq  c_{19}.
$$

\nind {\it Step 3: There exists $c_{20}>0$ such that
 \begin{equation}\label{15-09-3}
   \left|\int_{A}\frac{I_z(y)}{|y|^{1+2s}}dy\right|\leq  c_{20},
 \end{equation}
where $A=(-\infty,-\frac32)\cup (\frac32,+\infty)$.}
In fact, for $y\in A$, we observe that there exists $c_{21}>0$ such that
$I_z(y)\leq  c_{21}$ and
$$
   \int_{A}\frac{I_z(y)}{|y|^{1+2s}}dy\leq  2\int_{\frac32}^{+\infty}\frac{c_{21}}{|y|^{1+2s}}dy\leq  c_{22}
$$
for some $c_{22}>0$.
Consequently, by (\ref{15-09-0})-(\ref{15-09-3}), there exists $c_{23}>0$ such that
$$
 (-\Gd)_{_{\BBR}}^s \tilde w(z)\geq  -\frac{c_{23}}{1+z^{1+2s}} \qquad \forall\, z\geq  2.
$$
Since $1-\tilde w\chi_{_{\R_-}}=1$ in $\R_+$ and $1-\tilde w\chi_{_{\R_-}}\leq 1$ in $\R_-$,
 we have also
$$(-\Gd)_{_{\BBR}}^s(1-\tilde w\chi_{_{\R_-}})(z)\geq 0 \qquad \forall\, z>0.$$
Therefore, we obtain that
\begin{equation}\label{15-09-4}
(-\Gd)_{_{\BBR}}^s  w(z)\geq  -\frac{c_{23}}{1+z^{1+2s}} \qquad \forall\, z\geq  2.
\end{equation}
Combining (\ref{14-09-5}) and (\ref{15-09-4}), we infer that there exists  $R_1\geq  R_0+2$ such that for $z>R_1$,
$$
\begin{array} {ll}
\displaystyle
   (-\Gd)_{_{\BBR}}^s w(z)-\frac1{2s}w'(z)z-\frac{1+\beta}{p-1}w(z)
   \geq \sigma_0w(z)-\frac{c_{23}}{1+z^{1+2s}}
  \\[4mm]\displaystyle\phantom{   (-\Gd)_{_{\BBR}}^s w(z)-\frac1{2s}w'(z)z-\frac{1+\beta}{p-1}w(z)}
   =w(z)\left(\sigma_0-\frac{c_{23}}{\ln(e+z^2)}\right)\
\\[4mm]\displaystyle\phantom{  (-\Gd)_{_{\BBR}}^s w(z)-\frac1{2s}w'(z)z-\frac{1+\beta}{p-1}w(z)}
  \geq 0.
\end{array}
$$
For $z\leq  R_1$,  there exists $c_{24}>0$ such that
$$(-\Gd)_{_{\BBR}}^s w(z)-\frac1{2s}w'(z)z-\frac{1+\beta}{p-1}w(z)\geq  -c_{24},$$
and
there exists $c_{25}>0$ dependent of $R_1$ such that
$$w(z)\geq  c_{25}.$$
Therefore, one can find $\Lambda_0>0$ such that for $\lambda\geq \Lambda_0$,
\begin{equation}\label{15-09-5}
  (-\Gd)_{_{\BBR}}^ w(z)-\frac1{2s}w'(z)z-\frac{1+\beta}{p-1}w(z)+\lambda^{p-1}w^p(z)\geq 0 \qquad \forall\, z\in\R,
\end{equation}
which, together with (\ref{15-09-00}), implies that (\ref{14-09-0}) holds true. This ends the proof.
\hfill$\Box$\medskip

\subsection{Solutions with initial trace $(\CS,0)$ }

\blemma {ball} Assume that $N\geq 1$ and $p>1+\frac{2s(1+\gb)}{1+2s}$. Then for any $R>0$ there exists a positive function $u=u_{\infty,B_R}$ minimal among the solutions of (\ref{In-15}) in $Q_\infty$, which satisfy
\bel{Bn-3}\BA {lll}\displaystyle
\lim_{t\to 0}u(t,x)=\infty&\quad{\rm uniformly\; in}\quad \overline B_R\\[2mm]
\displaystyle
\lim_{t\to 0}u(t,x)=0&\quad{\rm uniformly\; in}\quad  B^c_{R'}\quad{\rm for\ any\ }  R'>R.
\EA\ee
Furthermore, the mapping $R\mapsto u_{\infty,B_R}$ is increasing.
\es
\Proof By scaling we can assume that $R=1$ and we fix $\gl\geq\gl_0$. We denote by ${\bf e}_1$ the point with coordinates $(1,0,...,0)$ in $\BBR^N$. The function
\bel{Bn-4}(t,x)\mapsto w_{e_1}(t,x_1,x')=\gl t^{-\frac{1+\gb}{p-1}}W(t^{-\frac{1}{2s}}(x_1-1)),
\ee
is a super solution of (\ref{In-15}) in $Q_\infty$, which satisfies
\bel{Bn-5}\BA {lll}\displaystyle
(i)\qquad&\lim_{t\to 0}w_{e_1}(t,x_1,x')=\infty&\quad{\rm uniformly\; in}\quad (-\infty,1]\ti\BBR^{N-1},\\[2mm]
\displaystyle
(ii)\qquad&\lim_{t\to 0}w_{e_1}(t,x)=0&\quad{\rm uniformly\; in}\quad [1+\ge,\infty)\ti\BBR^{N-1}.
\EA\ee
Since equation (\ref{In-15}) is invariant under rotations and translations, for any $a\in\prt B_1$ there exists a rotation $\CR_a$ with center $0$ such that $\CR_a(a)=e_1$. Therefore, the function $(t,x)\mapsto w_a(t,x):=w_{e_1}(t,\CR_a(x)) $ is a solution of  (\ref{In-15}) in $Q_\infty$ and it satisfies
\bel{Bn-6}\BA {lll}\displaystyle
(i)\qquad&\lim_{t\to 0}w_{a}(t,x)=\infty&\qquad{\rm uniformly\; in}\quad \{x\in\BBR^N:\langle x,a\rangle\leq 1\},\\[2mm]
\displaystyle
(ii)\qquad&\lim_{t\to 0}w_{e_1}(t,x)=0&\qquad{\rm uniformly\; in}\quad \{x\in\BBR^N:\langle x,a\rangle\geq 1+ge\}.
\EA\ee
For $k\in\BBN^*$, let $u_{k\chi_{_{B_1}}}$ be the solution of
\bel{Bn-7}\BA {lll}
\prt_tu+(-\Gd)^s u+t^\gb u^p=0&\quad{\rm in }\quad Q_\infty\\[1mm]
\phantom{\prt_tu+(-\Gd)^s u+}
\!u(0,.)=k\chi_{_{B_1}}&\quad{\rm in }\quad \BBR^N.
\EA\ee
Then the sequence $\{u_{k\chi_{_{B_1}}}\}_k$ is increasing. For any $a\in\prt B_1$, $u_{k\chi_{_{B_1}}}\leq w_a$, the following limit exists,
$$u_{\infty,B_1}=\lim_{k\to\infty}u_{k\chi_{_{B_1}}},$$
and there holds
$$u_{\infty,B_1}\leq \inf\left\{w_a: a\in\prt B_1\right\}.$$
This solution $u$ is clearly minimal by construction and the monotonicity of the mapping $R\mapsto u_{\infty,B_R}$ follows.\qeda
\medskip

\nind\Remark In the previous result, the ball $B_R$ can be replaced by any closed convex set with a non-empty interior. If $a\in\prt K$, let
$H_a$ be an affine separation hyperplane, with outer normal vector ${\bf n}_a$ and
$$H_a^+=\{x\in  \BBR^N:\langle x-a,{\bf n}_a\rangle>0\}\quad{\rm and}\quad H_a^-=\{x\in  \BBR^N:\langle x-a,{\bf n}_a\rangle<0\}.
$$
The supersolutions $w_a$ are expressed by
$$(t,x)\mapsto w_a(t,x)=\gl t^{-\frac{1+\gb}{p-1}}W(t^{-\frac{1}{2s}}\langle x-a,{\bf n}_a\rangle)
$$
and have initial trace $(0,\overline H_a^-)$. Then we construct the minimal solution $u=u_{\infty,K}$ of (\ref{In-15}) with initial trace $(0,K)$ such that
\bel{Bn-8}\BA {lll}\displaystyle
(i)\qquad&\lim_{t\to 0}u(t,x)=\infty&\quad{\rm uniformly\; in}\quad K,\\[2mm]
\displaystyle
(ii)\qquad&\lim_{t\to 0}u(t,x)=0&\quad{\rm uniformly\; in}\quad \{x\in K^c:\dist (x,K)\geq \ge\}\quad\forall \ge>0.
\EA\ee
Furthermore, the mapping $K\mapsto u_{\infty,K}$ is nondecreasing. \medskip


\bprop {CS} Assume that $N\geq 1$ and $p>1+\frac{2s(1+\gb)}{1+2s}$. Then for any closed set $\CS$ such that $\overline{{\rm int}(\CS)}=\CS$ there exists a positive function $u=u_{\infty,S}$ minimal among the solutions of (\ref{In-15}) in $Q_\infty$ which satisfy
\bel{Bn-9}\BA {lll}\displaystyle
(i)\quad&\lim_{t\to 0}u(t,x)=\infty&\quad{\rm locally\;uniformly\; in}\quad \CS,\\[2mm]
\displaystyle
(ii)\quad&\lim_{t\to 0}u(t,x)=0&\quad{\rm locally\; uniformly\; in}\quad\{x\in \CS^c:\dist(x,\CS)\geq \ge\}\quad\forall\, \ge>0.
\EA\ee
In particular $Tr(u_{\infty,S})=(\CS,0)$. Furthermore,
\bel{Bn-9-1}\BA {lll}\displaystyle
u_{\CS,\infty}(t,x)
   \leq c_{9}t^{-\frac{1+\beta}{p-1}}\myfrac{\ln \left(e+t^{-\frac{1}{s}}(\dist (x,\CS))^2\right)}{1+t^{-\frac{1+2s}{2s}}(\dist (x,\CS))^{1+2s}}\qquad\forall\, (t,x)\in Q_\infty.
\EA\ee
\es

\nind\Proof We first assume that $\CS$ is compact, hence precompact, and for any $\gd>0$ there exists a finite number of points ${\xi_j}\in\CS$, $1\leq j\leq n_\gd$ such that
$$\CS\subset \bigcup_{j=I}^{n_\gd}\overline B_\gd(\xi_j):=\CS_\gd.
$$
Clearly the mapping $\gd\mapsto n_\gd$ is nondecreasing, furthermore we can choose the points ${\xi_j}$ such that $\gd\mapsto \CS_\gd$ is decreasing for the order relation of inclusion between sets. Since
$p>1$, the function
\bel{Bn-10}\BA {lll}\displaystyle
w_{\CS_\gd}:=\sum_{j=1}^{n_\gd}u_{\infty,\overline B_\gd(\xi_j)},
\EA\ee
is a supersolution of (\ref{In-15}) in $Q_\infty$  and by \rlemma {ball} it satisfies
\bel{Bn-11}\BA {lll}\displaystyle
(i)\qquad&\lim_{t\to 0}w_{\CS_\gd}(t,x)=\infty&\quad{\rm uniformly\; in}\quad \CS_\gd,\\[2mm]
\displaystyle
(ii)\qquad&\lim_{t\to 0}w_{\CS_\gd}(t,x)=0&\quad{\rm  uniformly\; in}\quad \{x\in \CS_\gd^c:\dist(x,\CS_\gd)\geq \ge\}\quad\forall \ge>0.
\EA\ee
For $k\in \BBN^*$ let $u_{k\chi_{_\CS}}$ be the solution of (\ref{In-15}) in $Q_\infty$ with initial data $k\chi_{_\CS}$. It exists since $\CS$ has a non-empty interior, and it coincides with  the solution of (\ref{In-15}) in $Q_\infty$ with initial data $k\chi_{_{{\rm int}(\CS)}}$. Clearly there holds
$u_{k\chi_{_\CS}}\leq w_{\CS_\gd}$ and the sequence $\{u_{k\chi_{_\CS}}\}_k$ is increasing, then there exists
$$u_{\infty,\CS}=\lim_{k\to\infty}u_{k\chi_{_\CS}}.$$
It is a positive solution of (\ref{In-15}) in $Q_\infty$ which tends to infinity on $\CS$, by construction, and satisfies $u_{\infty,\CS}\leq w_{\CS_\gd}$. This implies in particular that for any
$\ge>0$,
$$\lim_{t\to 0}u_{\infty,\CS}(t,x)=0\quad{\rm  uniformly\; in}\quad \left\{x\in \CS_\gd^c:\dist(x,\CS_\gd)\geq \ge\right\}.
$$
Since this holds for any $\gd,\ge>0$, the second assertion in (\ref{Bn-9}) follows.\smallskip

\nind If $S$ is unbounded, for any $\gr>0$ large enough, $\CS^\gr:=\CS\cap \overline B_\gr$ is a nonempty compact set and $\CS^\gr=\overline{\rm int}(\CS^\gr)$. Hence there exists a solution $u_{\infty,\CS^\gr}$ of (\ref{In-15}) in $Q_\infty$ with initial trace $(0,\CS^\gr)$. By construction $\gr\mapsto u_{\infty,\CS^\gr}$ is nondecreasing and converges to a nonnegative solution $u_{\infty,\CS}$ of (\ref{In-15}) in $Q_\infty$. Let $a=(a_1,...,a_N)\in \CS^c$ and $\gt>0$ such that
$$Q^\gt_{a}=\{x=(x_1,...,x_N):\abs{x_j-a_j}\leq\gt\}\subset \CS^c.
$$
We put
$$W_j(t,x_j)=\gl t^{-\frac{1+\gb}{p-1}}\left(W(t^{-\frac{1}{2s}}(x_j-a_j+\gt)+W(t^{-\frac{1}{2s}}(a_j+\gt-x_j)\right)
$$
with $\gl\geq\gl_0$, then $W_j$ is a supersolution of (\ref{In-15}) in $\BBR_+\ti\BBR$ and it satisfies
$$\BA {lll}\displaystyle
(i)\qquad&\displaystyle\lim_{t\to 0}W_j(t,x)=0\qquad&{\rm locally\; uniformly\; in}\quad (a_j-\gt,a_j+\gt),\\\displaystyle
(ii)\qquad&\displaystyle\lim_{t\to 0}W_j(t,x)=\infty\qquad&{\rm  uniformly\; in}\quad (-\infty,a_j-\gt]\bigcup[a_j+\gt,\infty).
\EA$$
Hence $W_{Q^\gt_{a}}(t,x)=\sum_jW_j(t,x)$ is a supersolution of (\ref{In-15}) in $Q_\infty$ and it satisfies
$$\BA {lll}\displaystyle
(i)\qquad\quad&\displaystyle\lim_{t\to 0}W_{Q^\gt_{a}}(t,x)=0\qquad&{\rm locally\; uniformly\; in}\quad Q^\gt_{a},\qquad\qquad\\\displaystyle
(ii)\qquad&\displaystyle\lim_{t\to 0}W_{Q^\gt_{a}}(t,x)=\infty\qquad&{\rm  uniformly\; in}\quad \BBR^N\setminus Q^{\gt}_{a}.\qquad\qquad
\EA$$
By construction $u_{\infty,\CS^\gr}\leq W_{Q^\gt_{a}}$ which implies $u_{\infty,\CS}\leq W_{Q^\gt_{a}}$. Hence $u_{\infty,\CS}$ satisfies (\ref{Bn-9}). The estimate from above can be made more precise (it does not depend from the fact that $\CS=\overline{\rm int}\,\CS$) using (\ref{bn-1}) since
\bel{Bn-11-1}W_{Q^\gt_{a}}(a)\leq 2N\gl t^{-\frac{1+\gb}{p-1}}\myfrac{\ln \left(e+ t^{-\frac{1}{s}}\gt^2\right)}{1+t^{-\frac{1+2s}{2s}}\gt^{1+2s}}.
\ee
If we take $\gt=\frac{\dist(a,\CS)}{\sqrt N}$, we obtain (\ref{Bn-9-1}).
Furthermore $u_{\infty,\CS}$ is clearly minimal as the limit of an increasing sequence of solutions with bounded initial data having compact support.
\qeda

\subsection{Proof of Theorem I  }

If $K\subset\CR$ is compact, then $\gn_K=\chi_{_K}\gn\in \mathfrak M_+^b(\CR)$; we extend it by zero and still denote by $\gn_K\in \mathfrak M^b(\BBR^N)$ its extension. For $\gr>0$, $\CS^\gr:=\CS\cap \overline B_\gr$ and for $\ell\in\BBN^*$, $\ell\chi_{_{\CS_\gr}}dx$ is a good measure. Since $\gn_K$ is a good measure, $\gn_K+\ell\chi_{_{\CS_\gr}}dx$ is a good measure by \rprop{add}. Then there exists a solution $u:=u_{\gn_K+\ell\chi_{_{\CS_\gr}}dx}$ of (\ref{g-1}) in $Q_\infty$ with initial data $\gn_K+\ell\chi_{_{\CS_\gr}}dx$ and it satisfies
\bel{Bn-12}\BA {lll}\displaystyle
\sup\left\{u_{\gn_K},u_{\ell\chi_{_{\CS_\gr}}dx}\right\}\leq u_{\gn_K+\ell\chi_{_{\CS_\gr}}dx}\leq u_{\gn_K}+u_{\ell\chi_{_{\CS_\gr}}dx}\leq
u_{\gn_K}+u_{\infty,\CS}.
\EA\ee
Since $(\ell,\gr)\mapsto u_{\ell\chi_{_{\CS_\gr}}dx}$ is increasing, we can let $\ell$ and $\gr$ go to infinity succesively and obtain that $u_{\gn_K+\ell\chi_{_{\CS_\gr}}dx}$ converges to a positive solution $\tilde u_K$ of (\ref{In-15}) in $Q_\infty$ and that
\bel{Bn-12'}\BA {lll}\displaystyle
\sup\left\{u_{\gn_K},u_{\infty,\CS}\right\}\leq \tilde u_K\leq
u_{\gn_K}+u_{\infty,\CS}.
\EA\ee
This estimate implies that $Tr(\tilde u_K)=(\CS,\gn_K)$. To end the proof we consider an increasing sequence $\{K_n\}$ of compact sets such that
$\bigcup_nK_n=\CR$. Then estimate $(\ref{Bn-12'})$ holds with $K$ replaced by $K_n$. Furthermore the sequences $\{u_{\gn_{K_n}}\}$ and
$\{\tilde u_{K_n}\}$ are increasing. In order to prove that the sequence $\{u_{\gn_{K_n}}\}$ converges to some solution $\tilde u_\gn$ of (\ref{In-15}) in $Q_\infty$ which admits $\gn$ as the regular part of its initial trace, for $R>0$ we write
$\gn_{K_n}=\chi_{_{\overline {B}_R}}\gn_{K_n}+\chi_{_{\overline {B_R}^c}}\gn_{K_n}$ (both solutions exist since $K_n$ is admissible). Then
\bel{Bn-13'}u_{\gn_{K_n}}\leq u_{\chi_{_{\overline {B}_R}}\gn_{K_n}}+u_{\chi_{_{ B_R^c}}\gn_{K_n}}\leq u_{\chi_{_{\overline {B}_R}}\gn}
+u_{\infty, B_R^c\cap K_n}\leq u_{\chi_{_{\overline {B}_R}}\gn}
+u_{\infty, B_R^c\cap \overline\CR},
\ee
which implies that the following limit exists and satisfies the upper estimate for any $R>0$,
\bel{Bn-13}
\displaystyle
\lim_{n\to\infty}u_{\gn_{K_n}}:=\tilde u_\gn\leq u_{\chi_{_{\overline {B}_R}}\gn_{K}}+u_{\infty, B_R^c\cap \overline\CR}.
\ee
In turns it implies
\bel{Bn-14}
\displaystyle
\sup\left\{\tilde u_{\gn},u_{\infty,\CS}\right\}\leq \lim_{n\to\infty}u_{{K_n}}:=\tilde u\leq \tilde u_{\gn}+u_{\infty,\CS}.
\ee
Furthermore, since $R>0$ in inequality $(\ref{Bn-13})$ we infer that $\gn$ is the regular part of the initial trace of $\tilde u_\gn$ (notice that the singular part is not empty since $\gn$ can be unbounded). Hence $Tr(\tilde u)=(\CS,\gn)$.
\qeda

\subsection{Proof of Corollary K, part (a)}

 If $\gn$ vanishes on Borel sets with zero ${\rm cap\,}^{\BBR^N}_{\frac{2s(1+\gb)}{p},p'}$-capacity, for any compact set $K\subset\BBR^N$, $\gn_K:=\chi_{_K}\gn$ vanishes also on the same Borel sets. Hence there exists a solution $u_{\gn_K}$ to $(\ref{g-1})$ with initial data $\gn_K$ (instead of $\gn$). Next we replace $K$ by an increasing sequence $\{\overline B_n\}_{n\in\BBN^*}$, and set
$\gn_n=\chi_{_{\overline B_n}}\gn$. Estimate $(\ref{Bn-13'})$ holds in the form
$$u_{\gn_{K_n}}\leq  u_{\chi_{_{\overline {B}_R}}\!\!\gn}
+u_{\infty, B_R^c}\qquad\forall\, n,R\geq 1.
$$
This implies that $\tilde u_\gn$ satisfies the same estimate for any $R>0$, which in turn implies that
$$\displaystyle\lim_{t\to\infty}\myint{\BBR^N}{}\tilde u_\gn(t,x)\gz(x) dx =\myint{\BBR^N}{}\gz d\gn\qquad\forall\, \gz\in C^2_0(\BBR^N).
$$
Hence $Tr(\tilde u_\gn)=(\{\emptyset\},\gn)$. The fact that $\tilde u_\gn\in L^1_{loc}(0,T;\CL^s(\BBR^N))$ follows from the upper estimate
$0\leq \tilde u_\gn\leq U_{p,\gb}$. \smallskip

Conversaly (and here we do not use the assumption $p>1+\frac{2s(1+\gb)}{1+2s}$), if $u\in L^1_{loc}(0,\infty;\CL^s(\BBR^N))$ is a solution with initial trace $Tr(u)=(\{\emptyset\},\gn)$, then $u\leq U_{p,\gb}$, by Theorem D.
We proceed as in the proof of \rprop{adm}. Let $K\subset B_R\subset \BBR^N$ be compact and $\Gth\in C^{\infty}_0(B_{2R})$ such that $0\leq\Gth\leq 1$ and $\Gth= 1$ on $\overline B_R$. Since
$$\norm{\Gth\phi}_{W^{k,p'}}\leq c(k,p)\norm{\phi}_{W^{k,p'}}\qquad\forall\, \phi\in C^{\infty}_0(R^N)
$$
for all $k\in\BBN^*$ by Leibnitz formula, it follows by interpolation that
\bel{Bn-15}\norm{\Gth\phi}_{L^{\frac{2s(1+\gb)}{p},p'}}\leq c(s,p)\norm{\phi}_{L^{\frac{2s(1+\gb)}{p},p'}}\qquad\forall\, \phi\in C^{\infty}_0(R^N).
\ee
If $T\in C^{\infty}(\BBR_+)$ satisfies
$$\displaystyle \sup_{t>0}\abs{t^{j-1}T(j)(t)}\leq L<\infty,\qquad\forall\, j=0,1,...,\ell:=\BA {lll}
\BBE\left[\frac{2s(1+\gb)}{p}\right]+1,\EA
$$
then, by {\it the smooth truncation theorem} (see \cite[Th. 3.3.3]{AdHe}),
\bel{Bn-16}\norm{T\left({\bf J}_{\frac{2s(1+\gb)}{p}}[\phi]\right)}_{L^{\frac{2s(1+\gb)}{p},p'}}\leq AL
\norm{{\bf J}_{\frac{2s(1+\gb}{p}}[\phi]}_{L^{\frac{2s(1+\gb)}{p},p'}}:=AL\norm{\phi}_{L^{\frac{2s(1+\gb)}{p},p'}}
\qquad\forall\, \phi\in C^{\infty}_0(R^N).
\ee
If we take in particular a function $T$ with value $1$ in $[1,\infty)$ we infer that if ${\rm cap\,}^{\BBR^N}_{\frac{2s(1+\gb)}{p},p'}(K)=0$, there exists a sequence $\{\gz_n\}\subset C^\ell_0(B_{2R})$ such that $0\leq\gz_n\leq 1$, $\gz_n=1$ on $K$ and
$\norm{\gz_n}_{L^{\frac{2s(1+\gb)}{p},p'}}\to 0$ as $n\to\infty$.
We set $\Gf_n=e^{-\gd t}\BBH_s[\gz_n]$ and take $\Phi_n^{p'}$ for test function. Then for any $\ge>0$ we have
$$\BA {lll}\myint{\BBR^N}{}(u\Phi_n^{p'})(1,x)dx+2p'\left(\myint{\ge}{1}\myint{\BBR^N}{}t^{-\frac{p'\gb}{p}}\left|(-\Gd)^s\Phi_n+\gd\Phi_n\right|^{p'}dx dt\right)^{\frac{1}{p'}}
\left(\myint{\ge}{1}\myint{\BBR^N}{}
t^\gb u^p\Phi_n^{p'}dxdt\right)^{\frac{1}{p}}\\[4mm]
\phantom{----))-----------------\ }
+\myint{\ge}{1}\myint{\BBR^N}{}
t^\gb u^p\Phi_n^{p'}dxdt\geq \myint{B_{2R}}{}(u\Phi_n^{p'})(\ge,x)dx.
\EA$$
When $\ge\to 0$, the right-hand side of the above inequality converges to $\myint{B_{2R}}{}\!\!\gz_n^{p'}d\gn\geq \gn(K)$.
Furthermore, we have that
$$\displaystyle\lim_{n\to\infty}\myint{\BBR^N}{}(u\Phi_n^{p'})(1,x)dx=0,
$$
by the dominated convergence theorem since $u\leq U_{p,\gb}$ and $\Phi_n^{p'}(1,x)\to 0$ for all $x\in\BBR^N$, and
$$\displaystyle\lim_{n\to\infty}\lim_{\ge\to 0}\myint{\ge}{1}\myint{\BBR^N}{}t^{-\frac{p'\gb}{p}}\left|(-\Gd)^s\Phi_n+\gd\Phi_n\right|^{p'}dx dt=0,
$$
as in the proof of \rprop{charac}. For the last term, we have
$$\myint{\ge}{1}\myint{\BBR^N}{}
t^\gb u^p\Phi_n^{p'}dxdt=\myint{\ge}{1}\myint{B_{2R}}{}
t^\gb u^p\Phi_n^{p'}dxdt+\myint{\ge}{1}\myint{B^c_{2R}}{}
t^\gb u^p\Phi_n^{p'}dxdt:=I_{\ge,n}+J_{\ge,n}.
$$
By assumption $t^\gb u^p\in L^1((0,1)\ti B_{2R})$, then $\displaystyle\lim_{n\to\infty}\lim_{\ge\to 0}I_{\ge,n}=0$ by
the dominated convergence theorem. Finally, since $u\leq U_{p,\gb}$ by Theorem H  and $H_s(t,x)\leq \myfrac{ct}{t^{1+\frac{N}{2s}}+|x|^{N+2s}}$ by (\ref{cks}), we obtain with various $c>0$ independent of $R$
$$\BA {lll}J_{\ge,n}\leq c\myint{\ge}{1}\myint{B^c_{3R}}{}t^{(1-p')\gb} \left(\myint{B_{2R}}{}\myfrac{dy}{t^{1+\frac{N}{2s}}+\abs{x-y}^{N+2s}}\right)^{p'}dx dt\\[4mm]
\phantom{J_{\ge,n}}
\leq cR^{Np'}\myint{\ge}{1}t^{(1-p')\gb}\myint{B^c_{3R}}{}\myfrac{dx}{\left(t^{1+\frac{N}{2s}}+\abs{x}^{N+2s}\right)^{p'}}
\\[4mm]
\phantom{J_{\ge,n}}
\leq cR^{Np'}\myint{\ge}{1}t^{(1-p')\gb+\frac{N}{2s}-\frac{(N+2s)p'}{2s}}
\myint{3Rt^{-\frac{1}{2s}}}{\infty}\myfrac{r^{N-1}}{\left(1+r^{N+2s}\right)^{p'}}\\[4mm]
\phantom{J_{\ge,n}}
\leq cR^{N-2s p'}\myint{\ge}{1}t^{(1-p')\gb}dt
\\[2mm]
\phantom{J_{\ge,n}}
\leq c\frac{R^{N-2s p'}}{(p'-1)\gb-1}.
\EA$$
(Note that the assumption $\gb<p-1$ is crucial). Hence $\displaystyle\lim_{n\to\infty}\lim_{\ge\to 0}I_{\ge,n}=0$, always by the dominated convergence theorem. This implies that $\gn(K)=0$.\qeda

\section{The subcritical case }
\setcounter{equation}{0}

For equation (\ref{In-15}), the subcritical case corresponds to the fact that
$$u_\infty(t,x)=V(t,x)=t^{-\frac{1+\gb}{p-1}}v(t^{-\frac{1}{2s}}x)\qquad\forall\, (t,x)\in Q_\infty,$$
where $v$  is the minimal positive solution of (\ref{In-19}).

\subsection{Proof of Theorem J}

\begin{proposition}\label{lm 4.1}
Assume that  $1+\frac{2s(1+\beta)}{N+2s}<p<1+\frac{2s(1+\beta)}N$ and $u$ is a nonnegative solution of  (\ref{In-22}) where $\CS\neq\emptyset$. Then
\begin{equation}\label{4.1}
u(t,x)\geq   \frac{c_{10}t^{-\frac{1+\beta}{p-1}}}{1+(t^{-\frac1{2s}}d(x,\mathcal{S}))^{N+2s}} \qquad \forall \, (t,x)\in Q_\infty.
\end{equation}
for some $c_{10}>0$.
\end{proposition}
{\bf Proof.}
By Theorem C, for any $x_0\in \mathcal{S}$,
$$u(t,x)\geq  u_\infty(t,x-x_0) =t^{-\frac{1+\gb}{p-1}}v(t^{-\frac{1}{2s}}(x-x_0))\qquad \forall \, (t,x)\in Q_\infty,$$
which implies that
\bel{sc1}u(t,x)\geq  t^{-\frac{1+\gb}{p-1}}\sup_{x_0\in\mathcal{S}}v(t^{-\frac{1}{2s}}(x-x_0)) \qquad \forall\,  (t,x)\in Q_\infty.\ee
The maximum of $V$ is achieved at $0$, hence, for any $x\in\CS$,
\bel{sc2}u(t,x)\geq  t^{-\frac{1+\gb}{p-1}}V(0)=c_{11}t^{-\frac{1+\gb}{p-1}}.\ee
If $x\in\mathcal{S}^c$, there exists $\bar x\in \mathcal{S}$ such that $\dist(x,\mathcal{S})=|x-\bar x|$. It follows from  \cite[Theorem 1.2]{CVW} that,
\bel{sc3}
  u(t,x)
   \geq \frac{c_{10}t^{-\frac{1+\beta}{p-1}}}{1+(t^{-\frac1{2s}}\dist(x,\mathcal{S}))^{N+2s}}.
\ee
Then (\ref{4.1}) holds true.   \hfill$\Box$\medskip

The next result shows that any closed set can be the singular set of the initial trace of a positive solution of (\ref{In-15}).

\bprop {Prop42} Assume that $1+\frac{2s(1+\beta)}{1+2s}<p<1+\frac{2s(1+\beta)}N$ and $\CS\subset \BBR^N$ is a nonempty closed set. Then
there exists a minimal solution $u:=u_{\CS,\infty}$ with initial trace $(\CS,0)$. Furthermore it satisfies (\ref{Bn-9-1}).
\es
\Proof We first notice that the condition $1+\frac{2s(1+\beta)}{1+2s}<p<1+\frac{2s(1+\beta)}N$ is equivalent to the conditions stated in Theorem J, i.e.
\bel{4.1-1}\BA {lll}
&(i)\qquad\qquad{\rm either}\ \ N=1\; {\rm and}\; 1+\frac{2s(1+\beta)}{1+2s}<p<1+2s(1+\beta),\qquad\qquad\\[2mm]
&(ii) \qquad\qquad{\rm or}\quad \ \  N=2,\  \frac 12\leq s<1\; {\rm and}\; 1+\frac{2s(1+\beta)}{1+2s}<p<1+s(1+\beta).\qquad\qquad
\EA\ee
Let $A:=\{z_n\}_{n\in\BBN}\}$ be a countable dense subset of $\CS$. For $k\in\BBN_*$, set
\bel{4.1-2}\gm_k=k\sum_{j=1}^k\gd_{z_j},
\ee
and let $u=u_{\gm_k}$ be the solution of
\bel{sc4}\BA{lll}
\prt_tu+(-\Gd)^s u+t^\gb u^p=0\quad&{\rm in }\quad Q_\infty\\[1mm]
\phantom{\prt_tu+(-\Gd)^s u,,}
\displaystyle
\!u(0,.)=\gm_k\quad&{\rm in }\quad \BBR^N.
\EA\ee
The sequence $\{u_{\gm_k}\}$ is increasing. If $a\in \CS^c$ and $d_a=\dist (a,\CS)$. By construction there holds
\bel{sc5}u_{\gm_k}\leq u_{B^c_{d_a}(a),\infty}.
\ee
Hence $u_{\gm_k}$ converges to some solution $\tilde u$ of (\ref{In-15}) in $Q_\infty$ which has zero initial trace on $B_{d_a}(a)$, for any $a\in \CS^c$ since (\ref{sc5}) still holds with $\tilde u$ instead of $u_{\gm_k}$, and  satisfies  $\tilde u\geq u_{z_j,\infty}$ for any $z_j\in A$. Hence $Tr(\tilde u)=(\CS,0)$. Estimate (\ref{Bn-9-1}) is independent of the geometry of $\CS$.\qeda
\medskip

\nind{\bf Proof of Theorem J}. It is similar to the one of Theorem I . We consider an increasing sequence of compact sets $\{K_k\}_{k\in\BBN^*}$ included in $\CR$ such that $\bigcup_k K_k=\CR$, set $\gn_k=\chi_{_{K_k}}\gn$ and $\tilde \gn_k=\gn_k+\gm_k$, where $\gm_k$  is defined by (\ref{4.1-2}). Then the solution of (\ref{In-15}) with initial data $\tilde \gn_k$ satisfies
\bel{sc6}
\sup\{u_{\gn_k}, u_{\gm_k}\}\leq u_{\tilde\gn_k}\leq u_{\gn_k}+u_{\gm_k}.
\ee
By the same argument as in the proof of Theorem I ,  the sequence $\{u_{\gn_k}\}$ is increasing and converges to a solution $u_\gn$ (\ref{In-15}) with initial trace $(\{\emptyset\},\gn)$. Hence the sequence $\{u_{\tilde\gn_k}\}$ which is also increasing.
converges to some solution $\tilde u$ of (\ref{In-15}) which satisfies
\bel{sc7}
\sup\{u_{\gn}, u_{\infty,\CS}\}\leq \tilde u\leq u_{\gn}+u_{\infty,\CS}.
\ee
Then $\tilde u$ has initial trace $(\CS,\gn)$.\qeda\medskip

The proof of Corollary K, part (b) is straightforward.\medskip

\nind\Remark We conjecture that the following more general version of Theorem J holds: {\it For any integer $\gk \in [1,N]$ any $p>1$ such that $1+\frac{2s(1+\gb)}{\gk+2s}<p<1+\frac{2s(1+\gb)}{N}$, any closed set $\CS$ contained in an affine plane of codimension $\gk$ and any bounded measure in $\CS^c$, there exists a solution $u$ of problem (\ref{In-22}).}  We notice that the condition on $p$ can be fulfilled for some  $p$ if and only if $N-\gk<2s$, hence either $\gk=N$ i.e. $\CS$ is a single point and no condition on $s$, or $\gk=N-1$ hence $\CS$ is contained in a straight line and $\frac12<s<1$.

\subsection{Proof of Theorem L}

The proof uses the method developed in \cite{ShVe}. The function
\bel{K1}
\tilde \phi(x)=\inf\{\phi(y):|y|\geq|x|\},
\ee
is radial, nondecreasing, smaller that $\phi$ and we write it as  $\tilde \phi(|x|)$. Furthermore it satisfies
\bel{K2}
\lim_{\abs x\to\infty}|x|^{-\frac{2s}{p-1}}\tilde \phi(x)>0.
\ee
We set
$$\tilde\phi_n(x)=\left\{\BA {lll}\tilde\phi(x)\qquad{\rm if}\;|x|\leq n\\[2mm]
\tilde\phi(n)\qquad{\rm if}\;|x|> n.
\EA\right.
$$
The existence of a solution $u_{\tilde\phi_n}$ of (\ref{In-15}) with initial trace $(\{\emptyset,\tilde\phi_n\}$ follows from the fact that $\BBH_s[\tilde\phi_n]$ exists by \cite{BSV} and that $\BBH_s[\tilde\phi_n]\geq u_{\tilde\phi_n\chi_{_{B_k}}}$ for any
$k\in\BBN^*$. Hence $u_{\tilde\phi_n}$ is the increasing limit of $u_{\tilde\phi_n\chi_{_{B_k}}}$ when $k\to\infty$.
it is obtained by replacing $\tilde\phi_n$  by $\tilde\phi_n\chi_{_{B_k}}$ and by letting $k\to\infty$.
The solution $Y_n$ of the differential equation $Y'(t)+t^\gb Y^p(t)=0$ with initial data $Y(0)=\tilde \phi(n)$ is expressed by
\bel{K4}
Y_n(t)=\myfrac{\tilde \phi(n)}{\left(1+\frac{p-1}{\gb+1}t^{\gb+1}(\tilde \phi(n))^{p-1}\right)^{\frac{1}{p-1}}}.
\ee
It larger than $u_{\tilde\phi_n}$. Let us denote by $w_n$ the solution of
\bel{K5}\BA {lll}
\prt_tw_n+(-\Gd)^{s}w_n+t^{\gb}Y^{p-1}_n(t)w_n=0\qquad&{\rm in}\,\,&(0,\infty)\ti\BBR^N\\[1mm]
\phantom{\prt_t+(-\Gd)^{,\!s}w_n+t^{\gb}Y^{p-1}}
w_n(0,x)=\tilde\phi_n(x)\qquad&{\rm in}\,\,&\BBR^N.
\EA\ee
Since $Y^{p-1}_n\geq u^{p-1}_{\tilde\phi_n}$, $w_n$ is smaller that $u_{\tilde\phi_n}$, moreover $w_n$ can be explicitely computed
\bel{K6}
w_n(t,x)=e^{-\int_0^ts^{\gb}Y^{p-1}_n(s)ds}\myint{\BBR^N}{}H_{s}(t,x-y)\tilde\phi_n(y) dy.
\ee
In particular,
$$\BA {lll}w_n(t,0)=e^{-\int_0^ts^{\gb}Y^{p-1}_n(s)ds}\left(\myint{B_n}{}H_{s}(t,y)\tilde\phi_n(y) dy+\myint{B^c_n}{}H_{s}(t,y)\tilde\phi_n(y) dy\right)\\[4mm]
\phantom{w_n(t,0)}
\geq cte^{-\int_0^ts^{\gb}Y^{p-1}_n(s)ds}\tilde \phi(n)\myint{B^c_n}{}\myfrac{dy}{t^{\frac{N+2s}{2s}}+|y|^{N+2s}}\\[4mm]
\phantom{w_n(t,0)}
\geq cte^{-\int_0^ts^{\gb}Y^{p-1}_n(s)ds}\tilde \phi(n)\myint{n}{\infty}\myfrac{r^{N-1}dr}{t^{\frac{N+2s}{2s}}+r^{N+2s}}\\[4mm]
\phantom{w_n(t,0)}
\geq ctn^{-2s}e^{-\int_0^ts^{\gb}Y^{p-1}_n(s)ds}\tilde \phi(n).
\EA$$
Next
$$\BA {lll}-\myint{0}{t}s^{\gb}Y^{p-1}_n(s)ds=-\myint{0}{t}\myfrac{s^{\gb}(\tilde \phi(n))^{p-1}ds}{1+\frac{p-1}{\gb+1}t^{\gb+1}(\tilde \phi(n))^{p-1}}=-\myfrac{1}{p-1}\ln\left(1+\frac{p-1}{\gb+1}t^{\gb+1}(\tilde \phi(n))^{p-1}\right).
\EA$$
We write
$$\BA {lll}
tn^{-2s}e^{-\int_0^ts^{\gb}Y^{p-1}_n(s)ds}\tilde \phi(n)=
e^{\ln\left(t\tilde\phi(n)\right)-2s\ln n-\frac{1}{p-1}\ln\left(1+\frac{p-1}{\gb+1}t^{\gb+1}(\tilde \phi(n))^{p-1}\right)}\\[2mm]
\phantom{tn^{-2s}e^{-\int_0^ts^{\gb}Y^{p-1}_n(s)ds}\tilde \phi(n)}
=e^{\frac{1}{p-1}\ln\left(\frac{t^{p-1}(\tilde \phi(n))^{p-1}}{n^{2s (p-1)}\left(1+\frac{p-1}{\gb+1}t^{\gb+1}(\tilde \phi(n))^{p-1}\right)}\right)}.
\EA$$
If we take $t=t_n$ such that
\bel{K7}\displaystyle\lim_{n\to\infty}t_n^{\frac{\gb+1}{p-1}}\tilde \phi(n)=+\infty,
\ee
the expansion of the term in the logarithm gives
\bel{K8}\frac{t^{p-1}(\tilde \phi(n))^{p-1}}{n^{2s (p-1)}\left(1+\frac{p-1}{\gb+1}t^{\gb+1}(\tilde \phi(n))^{p-1}\right)}
=\myfrac{\gb+1}{p-1}n^{-2s(p-1)}t_n^{p-2-\gb}(1+o(1))\quad{\rm as }\;n\to\infty.
\ee
If besides (\ref{K7}) it is assumed that
\bel{K9}\displaystyle\lim_{n\to\infty}n^{-2s(p-1)}t_n^{p-2-\gb}=+\infty,
\ee
we infer that $w_n(0,t_n)\to\infty$ as $n\to\infty$. Clearly the origin can be replaced by any $z\in\BBR^N$ and the previous calculation shows that this limit is uniform for $z$ belonging to compact sets on $\BBR^N$. Since $u_{\tilde\phi_n}\geq w_n$, we infer that
 \bel{K10}\displaystyle\lim_{n\to\infty}u_{\tilde\phi_n}(t_n,z)=+\infty\Longrightarrow
 \lim_{n\to\infty}u_{\tilde\phi_n}=u_\infty=U_{p,\gb},
\ee
by using (\ref{In-17}).
\qeda

\setcounter{equation}{0}
\section{Appendix: symmetry and monotonicity results}

The following is a variant of the maximum principle which will be used in the sequel.
\begin{lemma}\label{MP}
Let $R,T>0$, $\gd\in [0,T)$ and $Q$ be a  domain of $Q_\infty$ such that $\overline Q\subset (\gd,T)\times B_R$. Assume that $h\geq 0$ in $Q$
and $\psi\in C(\bar Q)$ satisfies
\begin{equation}\label{mp1}
\arraycolsep=1pt
\begin{array}{lll}
\partial_t\psi+(-\Delta)^s \psi + h(t,x)\psi\geq 0\quad  {\rm in}\ \ Q\\[2mm]
 \phantom{ \psi_t+(-\Delta)^s \psi + h(t,x)\psi }
\psi\geq 0\quad  {\rm in}\ \ ([\gd,T)\times B_R)\setminus Q.\\[2mm]
\end{array}
\end{equation}
Then $\psi$ is nonnegative in $[\gd,T)\times B_R$.
\end{lemma}
{\bf Proof.} Let $\ge\in (0,T-\gd]$. We first claim that $\psi$ is nonnegative in $[\gd,T-\ge]\times B_R$. If it does not hold,
and since $\psi\geq 0$ in $([\gd,T)\times B_R)\setminus Q$, then
there exists $(t_0,x_0)\in Q\cap([\gd,T-\ge]\times B_R)$ such that
$$\psi(t_0,x_0)=\min_{(t,x)\in[\gd,T-\ge]\times B_R}\psi(t,x)<0.$$
Then $\partial_t\psi(t_0,x_0)\leq 0$ and
$(-\Delta)^s \psi(t_0,x_0)<0$. Since $h\geq 0$ in $Q$ and $(t_0,x_0)\in Q$, there holds
$$\partial_t\psi(t_0,x_0)+(-\Delta)^s\psi(t_0,x_0)+h(t_0,x_0) \psi(t_0,x_0)<0,$$
which is a contradiction. Thus, $\psi$ is nonnegative in $[\gd,T-\ge]\times B_R$. Since $\ge$ is arbitrary, the result follows. Notice that we can take $R=\infty$ in the above proof provided $Q$ is a bounded domain.$\phantom{--}$\qeda
\medskip

Next we prove the following result.

\bprop {sym}Let $N\geq 1$, $\gb>-1$, $p>1$ and $g\in C(\BBR^N)$ be a nonnegative continuous radially symmetric and nonincreasing function which tends to $0$ when $|x|\to\infty$. If $u\in L^1_{loc}(0,\infty;\CL^s(\BBR^N)\cap C(\overline Q_\infty)$  is a nonnegative solution of (\ref{In-15}) in $Q_\infty$ which converges to $g$ uniformly when $t\to 0$, then $u$ is radially symmetric and nonincreasing.
\es
\Proof Since $u\in L^1_{loc}\left(0,\infty;\CL^s(\BBR^N)\right)\cap C(\overline Q_\infty)$, it is bounded from above  by $\BBH_s[g]$ and uniqueness holds as for the linear equation \cite {BSV}. Since the initial data is radially symmetric and the equation is invariant by rotations in $\BBR^N$, $u(t,.)$ is also radially symmetric. Because of uniqueness and stability, it is sufficient to prove the result for a function $u$ which initial data is obtained from the previous one by multiplying it by a smooth, even, nonincreasing and nonnegative function with compact support.
The corresponding solution of (\ref{In-15}) in $Q_\infty$, still denoted by $u$, is smooth in $Q_\infty$ and bounded from above by $\BBH_s[g]$. Hence it satisfies
\begin{equation}\label{mp2}\BA{llll}
(i)\displaystyle\qquad\qquad &\lim_{t\to\infty}u(t,x)=0\qquad&{\rm uniformly\; in}\;\ x\in\BBR^N,\qquad\qquad\qquad\qquad   \\[2mm]
(ii)\displaystyle\qquad\qquad &\lim_{|x|\to\infty}u(t,x)=0\qquad&{\rm uniformly\; in}\;\ t\in\BBR_+,\qquad\qquad\qquad\qquad\\[2mm]
(iii)\displaystyle\qquad\qquad &\lim_{t\to0}u(t,x)=g(x)\qquad&{\rm uniformly\; in}\;\ x\in\BBR^N.\qquad\qquad\qquad\qquad
\EA\end{equation}
 \smallskip
Next we use a moving plane method (see \cite{PQR} for other applications). For $\lambda\in\R$, we set $x_\lambda=(2\lambda-x_1,x')$  if $x=(x_1,x')\in\R^N$,
\begin{equation}\label{03-10-0}
  \Sigma_\lambda=\{x=(x_1,x')\in \R^N\  |\  x_1<\lambda\}
\end{equation}
and
$$T_\lambda=\{x=(x_1,x')\in \R^{N}\ |\  x_1=\lambda\}.$$
We observe that if $\lambda>0$, then $\{x_\lambda\  |\ x\in\Sigma_\lambda\}=\{x\in \R^N  |\ x_1>\lambda\}$ and
\begin{equation}\label{01-10-1}
|x_\lambda|>|x|\quad \ {\rm for}\ x\in \Sigma_\lambda.
\end{equation}

\nind {\it We claim that for any $\lambda>0$},
\begin{equation}\label{d3}
u(t,x)\geq  u(t,x_\lambda) \qquad \forall\,  (t,x)\in \BBR^*_+\times\Sigma_\lambda.
\end{equation}
Set $\varphi(t,x)=u(t,x)-u(t,x_\lambda)$ and suppose that (\ref{d3}) does not hold. Because of (\ref{mp2})  there holds $\displaystyle\lim_{|x|\to\infty}\varphi(t,x)=0$   uniformly with respect to $t\geq 0$, $\displaystyle\lim_{t\to\infty}\varphi(t,x)=0$ uniformly with respect to $x\in\BBR^N$ and $\displaystyle\lim_{t\to 0}\varphi(t,x)=g(x)-g(x_\gl)\geq 0$ uniformly with respect to $x\in\BBR^N$. It follows that there exists
$\varepsilon_0>0$ and $(t_0,x_0)\in\BBR^*_+\times\Sigma_\lambda$ such that
\begin{equation}\label{eq2.4}
\varphi(t_0,x_0)=\min_{(t,x)\in\overline \Gs_\gl}\varphi(t,x)=-\varepsilon_0<0.
\end{equation}
The function $\phi$ satisfies
\begin{equation}\label{d4}\BA {lll}
\prt_t\gf+(-\Gd)^s\gf+ h(t,x)\gf=0\qquad{\rm in}\;\;Q_\infty,
\EA
\end{equation}
for some $h(t,x)\geq 0$, and it has initial data $\gf(0,x)=g(x)-g(x_\gl)$ in $\BBR^N$. Take $\ge\in (0,\ge_0)$ and set $\gf_\ge=\gf+\ge$. Using (\ref{mp2}) we see that there exists $T_0>t_0>0$ and $R_0>|x_0| >0$ such that $\gf_\ge(t,x)\geq 0$ for $(t,x)\in \left([T,\infty)\times\BBR^N\right)\bigcup \left([0,\infty)\times B_{R}^c\right)$, for all $T\geq T_0$ and $R\geq R_0$. Furthermore there exists $\gd_0\in (0,t_0)$ such that for any  $\gd\in (0,\gd_0)$ such that $\gf_\ge(t,x)\geq 0$ for $(t,x)\in [0,\gd)\times\left(\BBR^N\cap \Gs_\gl\right)$.
We set
$$Q=\Gs_\gl\cap (\gd,T_0)\times B_{R_0}.
$$
We apply Lemma \ref{MP} in $[\frac\gd 2, T)\times B_R$ and conclude that $\gf_\ge\geq 0$ in  $[\frac\gd 2, T)\times B_R$, which contradicts the fact that $\gf_\ge(t_0,x_0)=\ge-\ge_0<0$. Hence (\ref{d3}) holds. Since $\gl>0$ is arbitrary, this implies in particular by continuity that
\begin{equation}\label{d5}
\myfrac{\partial u}{\partial x_1}(t,x_1,x')\leq 0\qquad\forall\, (t,x_1,x')\in \BBR_+\times \BBR_+\times \BBR^{N-1}.
\end{equation}
Similarly, we can get that
\begin{equation}\label{d6}
\myfrac{\partial u}{\partial x_1}(t,x_1,x')\geq 0\qquad\forall\, (t,x_1,x')\in \BBR_+\times \BBR_-\times \BBR^{N-1}.
\end{equation}
Since $u(t,x)$ is radially symmetric with respect to $x$, it implies that $u(t,x)\geq u(t,x')$ if $ |x|\leq |x'| $, which ends the proof.
\qeda

\end{document}